\documentclass[11pt]{report}
\usepackage{amsmath,amssymb,amsbsy,amsfonts,amsthm,latexsym,amsopn,amstext,amsxtra,epic, euscript,amscd,indentfirst,verbatim}
\usepackage{hyperref}
\usepackage{xcolor}
\usepackage{graphicx}
\usepackage{float} % To control the location of figures

\newtheorem*{definition*}{Definition}

\newcommand{\comments}[1]{} 
\DeclareMathOperator{\Ric}{Ric}
\DeclareMathOperator{\Hess}{Hess}

\begin{document}

\title{An Illustrated Introduction to the Ricci Flow}

\author{Gabriel Khan}

%\maketitle

\begin{titlepage}
   \begin{center}
       \vspace*{1cm}

       \LARGE{An Illustrated Introduction to the Ricci Flow}
       
\vspace{.53in}
       \large{Gabriel Khan}
       
\vspace{.3in}
\large{\today}

       \vfill
            
    \begin{figure}[H]
    \centering
\includegraphics[width=.9\linewidth]{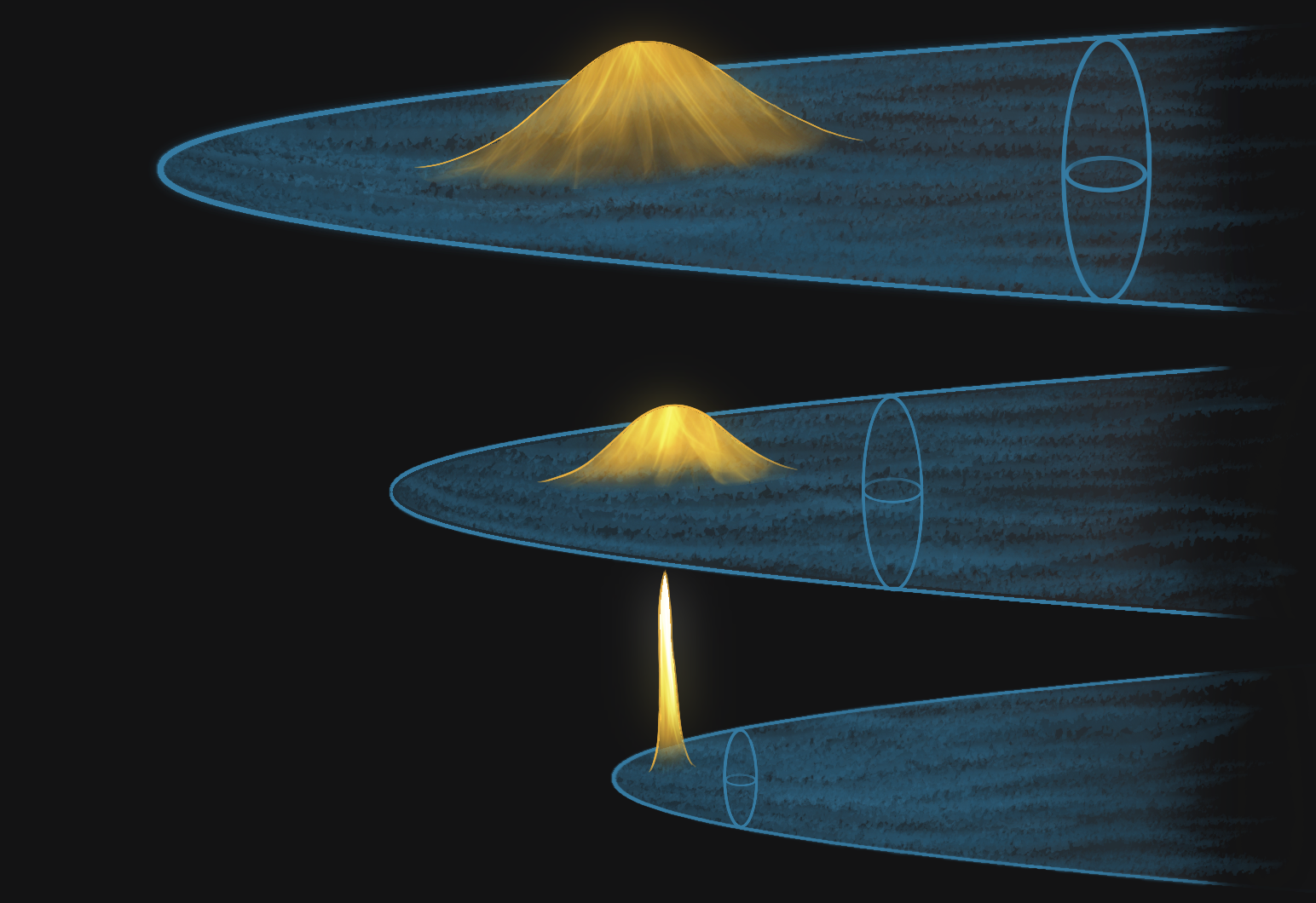}
    \caption{A Ricci flow with a conjugate heat flow\protect\footnotemark }
    \label{fig:Bryant soliton}
\end{figure}

\footnotetext{This image depicts a Bryant soliton with a backwards heat flow. It is adapted from Figure 6 of \cite{bamler2021recent}.}
            
   \end{center}
\end{titlepage}

\tableofcontents

\section*{Introduction}

The Ricci flow is one of the most important topics in differential geometry, and a central focus of modern geometric analysis.  In this paper, we give an illustrated introduction to the Ricci flow for a general audience. In particular, we do not assume any background in differential geometry or differential topology, only that the reader has taken multivariate calculus. For those who are familiar with geometric analysis, we have provided additional details in the footnotes. However, these notes are not necessary to follow the main discussion, so can be ignored on an initial reading.

%It is particularly famous for its role in Perelman's proofs of the Poincar\'e and Geometrization conjecture. 

%a historical account for its development and describe its behavior in the context of geometric heat flows. The goal is to provide an informal and illustrated introduction which is accessible to a general audience. In particular, we do not assume that the reader has any background of differential geometry or differential topology, only that you are familiar with multivariate calculus. 

%In the second half of the paper, we will provide an explanation for what the Ricci flow actually is. One very useful aphorism for the Ricci flow is that it is a ``heat equation for curvature." This statement is not entirely correct, but is a useful starting point and our goal in this paper is to explain how the Ricci flow behaves similarly to (and differently from) a heat equation. 

\chapter*{A short history of the Poincar\'e conjecture}
\addcontentsline{toc}{chapter}{A short history of the Poincar\'e conjecture}

\begin{figure}[H]
    \centering
\includegraphics[width=.9\linewidth]{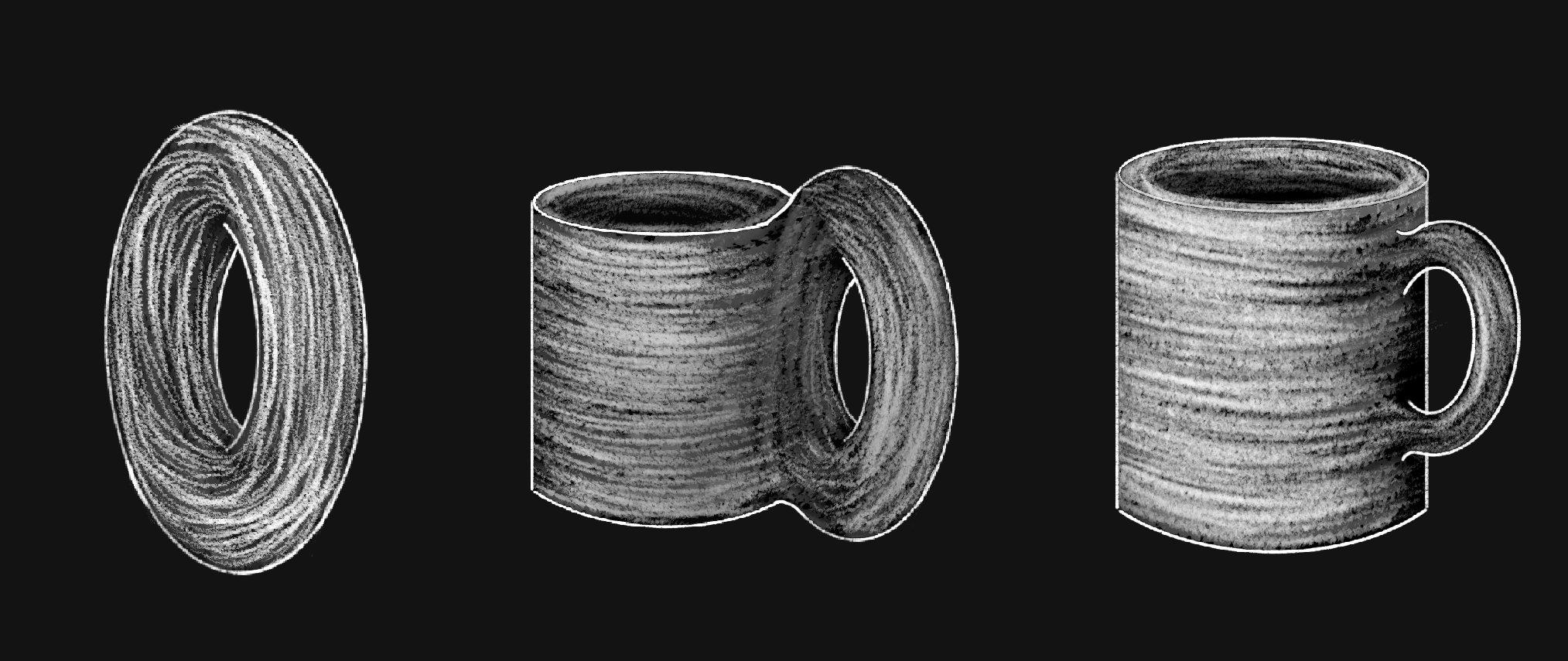}
\caption[Text]{A donut is topologically the same as a coffee mug}
\label{Donut deformation}
\end{figure}

The Ricci flow is most famous for its role in the proof of the Poincar\'e conjecture, so we start our discussion with a historical account of the Poincar\'e conjecture and some related problems in topology. 

In 1905, Henri Poincar{\'e} published the seminal paper \emph{Analysis Situs} \cite{poincare1895analysis}, which laid the foundations for what is known today as \emph{topology}.\footnote{Some topological ideas appeared  in Leonard Euler's earlier work on the Bridges of Konigsburg, but Analysis Situs is a foundational work in this area.} Roughly speaking, topology studies ``rubber-sheet geometry," where spaces are allowed to bend and stretch without tearing. Warping a space will change distances, areas, etc. so the properties studied in topology are quite different from those in traditional geometry. However, one central focus is to find features which do not change when a space is deformed continuously.

While studying this topic, Poincar\'e conjectured that the three-sphere is the only compact simply-connected three-dimensional manifold.\footnote{This is actually not Poincar\'e's original question, but a refinement after counter-examples were found to some earlier versions of the conjecture.} To explain the meaning of this conjecture, it is helpful to first consider the corresponding result for two-dimensional surfaces, which is a special case of the \emph{uniformization theorem}.

\section*{The uniformization theorem}

\begin{figure}[H]
    \centering
\includegraphics[width=.9\linewidth]{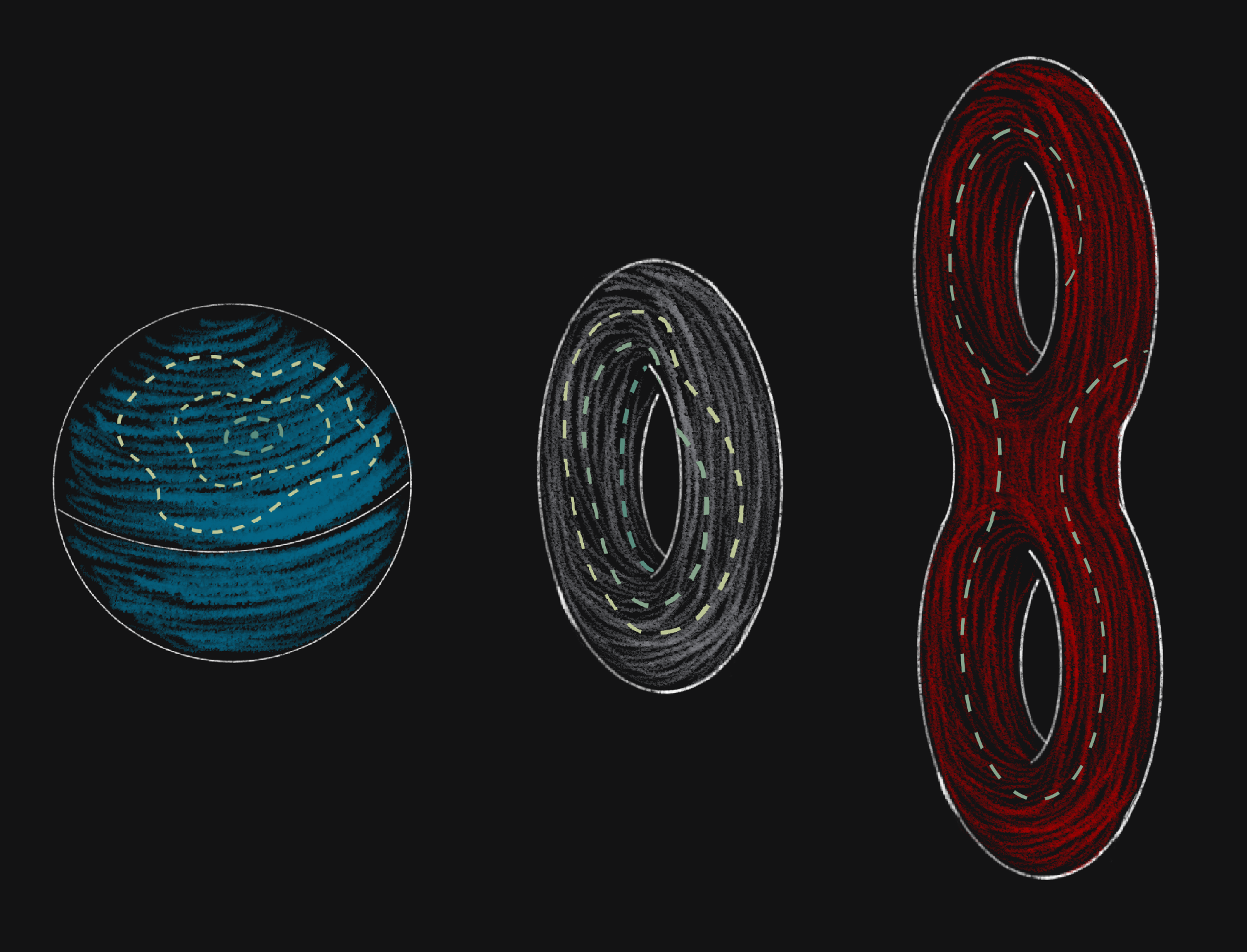}
\caption[Text]{Three surfaces and some loops on them\protect\footnotemark}
\label{Simply connected}
\end{figure}

\footnotetext{Throughout this paper we will use the color blue to indicate positive curvature, red for negative curvature and gray vanishing curvature.}

The uniformization theorem is a celebrated result which classifies the possible geometries on a two-dimensional surface \cite{poincare1908uniformisation}. It states any surface can be deformed into one of three special types of geometries, depending on the number of holes that it has.\footnote{The actual statement of the uniformization theorem is a bit stronger, and says that this deformation can be done in a conformal way. The theorem also gives a classification of non-orientable surfaces.}
Surfaces without any holes can be deformed into a round sphere. A surface with one hole (i.e., a donut) can be made into a flat space\footnote{Here, a surface is said to be \emph{flat} if its curvature vanishes and hyperbolic if its sectional curvature is identically $-1$. Later on, we will give a definition for curvature. } and any surface with more than one hole admits a hyperbolic geometry.

At first, it might be somewhat hard to visualize the latter two types of geometries. The reason for this is that when we think of surfaces, it is natural to visualize them as living in three-dimensional Euclidean space.\footnote{Here, by ``living inside three-dimensional space," we mean smoothly immersed in $\mathbb{R}^3$. For readers who are familiar with the differential geometry of curves and surfaces, it is a good exercise to prove that compact surfaces with non-positive curvature cannot be immersed in $\mathbb{R}^3$. As a hint, suppose the surface contains the origin in its interior. What can we say about the curvature of the point furthest from the origin?} However, there is no way to put a flat donut or a hyperbolic surface in $\mathbb{R}^3$. To get around this roadblock, we must consider the \emph{intrinsic} geometry of a surface, without assuming that it lies in some ambient Euclidean space. This can be a major conceptual challenge, but there is an intuitive way to visualize a flat donut using a classic arcade game.

\begin{figure}[H]
    \centering
\includegraphics[width=.9\linewidth]{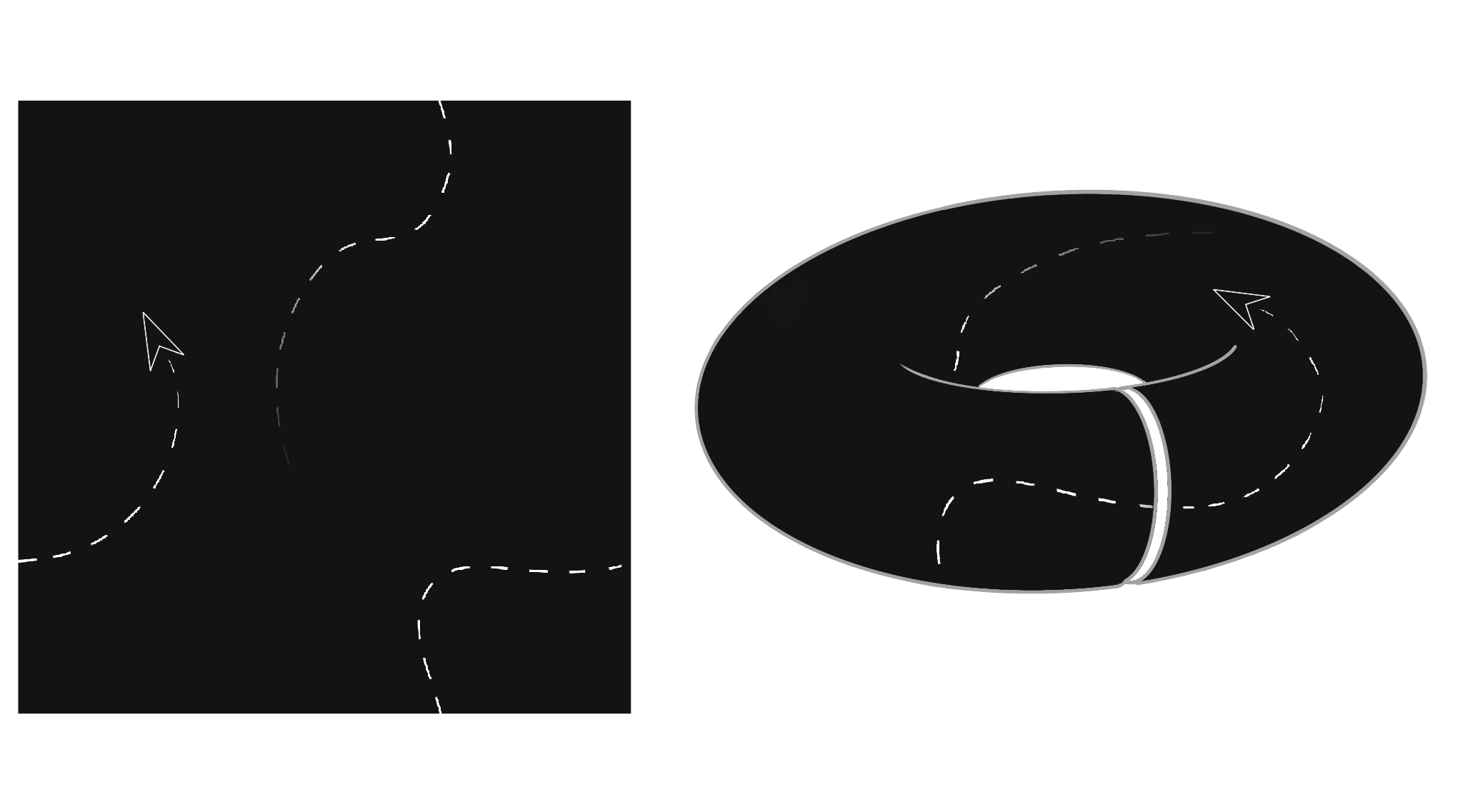}
\caption[Text]{Asteroids on a donut}
\label{Asteroids space}
\end{figure}

  In the game Asteroids, the player pilots a ship through an asteroid field on the screen, but there is a catch. When the ship flies across the right edge of the screen, it ends up on the left side of screen and when it flies across the top, it ends up at the bottom. It is a worthwhile exercise to show that a space like this is topologically equivalent to a standard donut.

It turns out that there are flat donuts in four-dimensional Euclidean space. Furthermore, there is a very deep theorem by John Nash which states that any possible geometry on a surface can be realized if we consider the surface as living in 51-dimensional Euclidean space \cite{nash1956imbedding}.\footnote{In general, any $n$-dimensional Riemannian manifold can be isometrically embedded in $\mathbb{R}^N$ where $N =n(n+1)(3n+11)/2$. The dimension $51$ comes from applying this formula to $n=2$. Mikhail Gromov proved that any compact surface can actually be isometrically embedded in $\mathbb{R}^5$, so for our purposes $46$ of these $51$ dimensions are redundant \cite{gromov1986partial}.} 

\section*{The Poincar\'e conjecture}

The uniformization conjecture has a straightforward, but very important, corollary relating the topology of a surface and the behavior of loops on it. As shown in Figure \ref{Simply connected}, if we consider a surface with a positive number of holes and draw a loop around the hole, it is not possible to shrink the loop to a point without cutting it or leaving the surface. 

As a result, one consequence of the uniformization theorem is that if we are given a surface where every loop drawn in it can be shrunk down continuously to a point, then the space is topologically equivalent to a sphere. Spaces where every loop can be contracted to a point are said to be \emph{simply-connected}, so another way to state this is that any simply-connected surface is topologically equivalent to a sphere.  The \emph{Poincar\'e conjecture} is the three-dimensional version of this statement. In other words, if we are given a closed three-dimensional space\footnote{Throughout the rest of the paper, we will use the word ``space" to mean ``compact smooth manifold" so as to not introduce unnecessary terminology.} which is bounded and simply-connected, then the space is topologically equivalent to a three-dimensional sphere.

\begin{figure}
    \centering
\includegraphics[width=.9\linewidth]{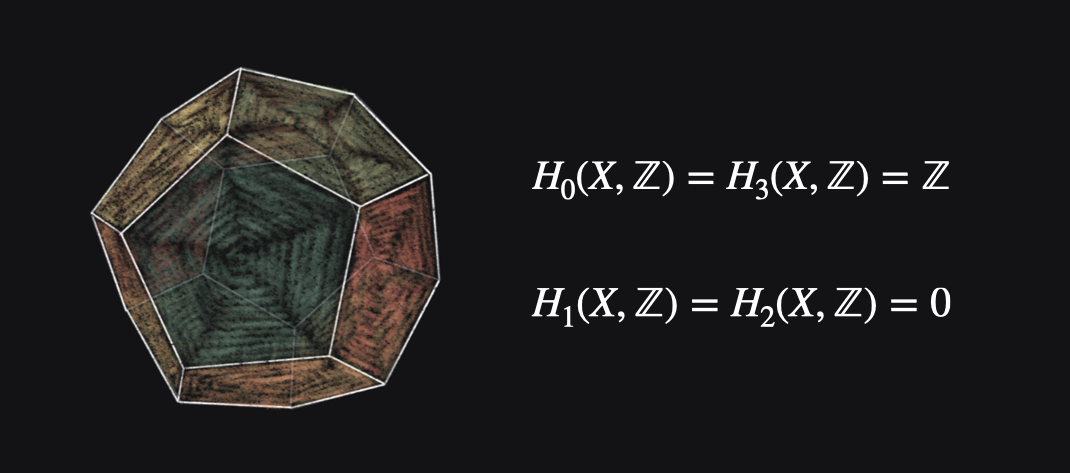}
\caption{A counter-example to an early version of the Poincar\'e conjecture\protect\footnotemark}
\label{Dodecahedron}
\end{figure}

\footnotetext{More precisely, if we take a solid dodecahedron and identify the opposite sides with a minimal twist, what results is a \emph{homology sphere}, which is a manifold whose homology groups are the same as the sphere but which is topologically distinct from a sphere.}

 This conjecture (and its higher dimensional analogues) motivated much of the early work in topology. The three-dimensional case in particular proved itself to be extremely subtle and intractable problem (see \cite{stallings2016not} for some idea of why this is the case). It attracted the attention of many leading mathematicians who developed various tools in their attempts to solve it. 
  As a result of all these efforts, by the 1980s the Poincar\'e conjecture was settled in all dimensions except for three.\footnote{There is an important subtlety here. There is a both a topological and a smooth version of the Poincar\'e conjecture, which considers whether spaces are continuously equivalent to a sphere or smoothly equivalent to the round sphere. The former is what is traditionally known as the Poincar\'e conjecture. Stephen Smale proved the topological Poincar\'e conjecture in dimensions greater than four \cite{smale1962structure} in 1962 and Michael Freedman proved the four-dimensional case \cite{freedman1982topology} in 1982. On the other hand, John Milnor showed that there are exotic smooth structures on the seven sphere, and thus the smooth Poincar\'e conjecture is false in dimension seven \cite{milnor1959differentiable}, which started a long line of work study the possible smooth structures on spheres. All three were awarded Fields medals for their respective work. At this point, there is one major question remaining, which is whether there exist exotic smooth structures on the four-dimensional sphere.} However, Poincar\'e's original conjecture appeared to be outside the reach of the standard tools of low-dimensional topology, and it seemed that a new approach would be needed to solve it.
  
  \section*{The Geometrization Conjecture}
  
  \begin{figure}[H]
    \centering
\includegraphics[width=.9\linewidth]{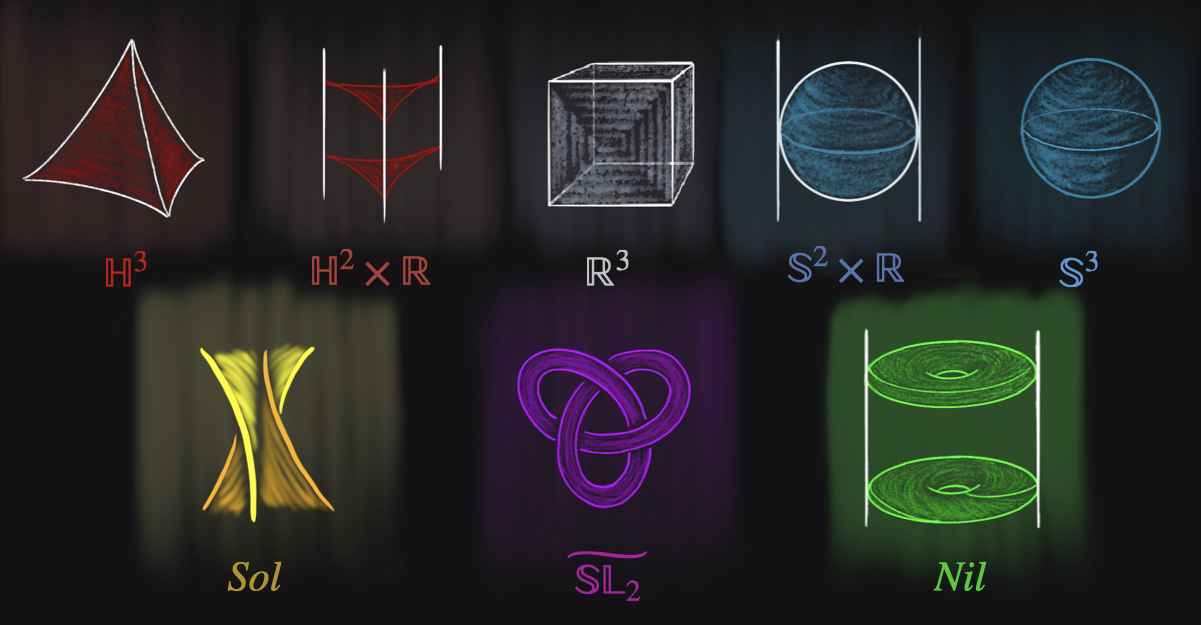}

    \caption[Thurston Geometries]{Thurston's Eight Geometries\protect\footnotemark}
    \label{Thurston's Geometries}
\end{figure}

  \footnotetext{In this illustration, the top row illustrates the five geometries which are either spaces of constant curvature or products thereof. The illustrations of the bottom row deserve further explanation.		
\begin{enumerate}
    \item The depiction of Sol illustrates the metric
   \[g= e^z dx^2 + e^{-z} dy^2+dz^2,\]
  which expands exponentially in the $x$-direction and contracts exponentially in the $y$-direction as the $z$-coordinate increases.
  \item The fundamental example of a space which has $ \widetilde{ SL}(2,\mathbb{R})$-geometry is the complement of the trefoil knot in the three-sphere.
  \item
Nil has the geometry of a torus bundle over a circle where the monodromy is given by a Dehn twist \cite{coulon2020non}. Here, the vertical lines depict the universal cover of a circle so that the two tori appear at different heights.
\end{enumerate} 
   }
  
In the mid-1970s, William Thurston began a line of work to understand the geometry of three-dimensional spaces, with an emphasis on those with hyperbolic geometry. Building from this work, he proposed a classification for the possible geometries of three-dimensional spaces similar to how the uniformization theorem classifies the geometry of two-dimensional surfaces. More precisely, he found eight canonical three-dimensional geometries and conjectured that given any three-dimensional space, there was a natural way to cut it into pieces which admit one of the geometries. Thurston was able to prove this conjecture for a fairly broad class of spaces and was awarded the Fields Medal in 1982 for this work. His conjecture, known as Thurston's Geometrization conjecture, implied the Poincar\'e conjecture and gave a new avenue to attack Poincar\'e's question.

\begin{figure}[H]
    \centering
\includegraphics[width=.9\linewidth]{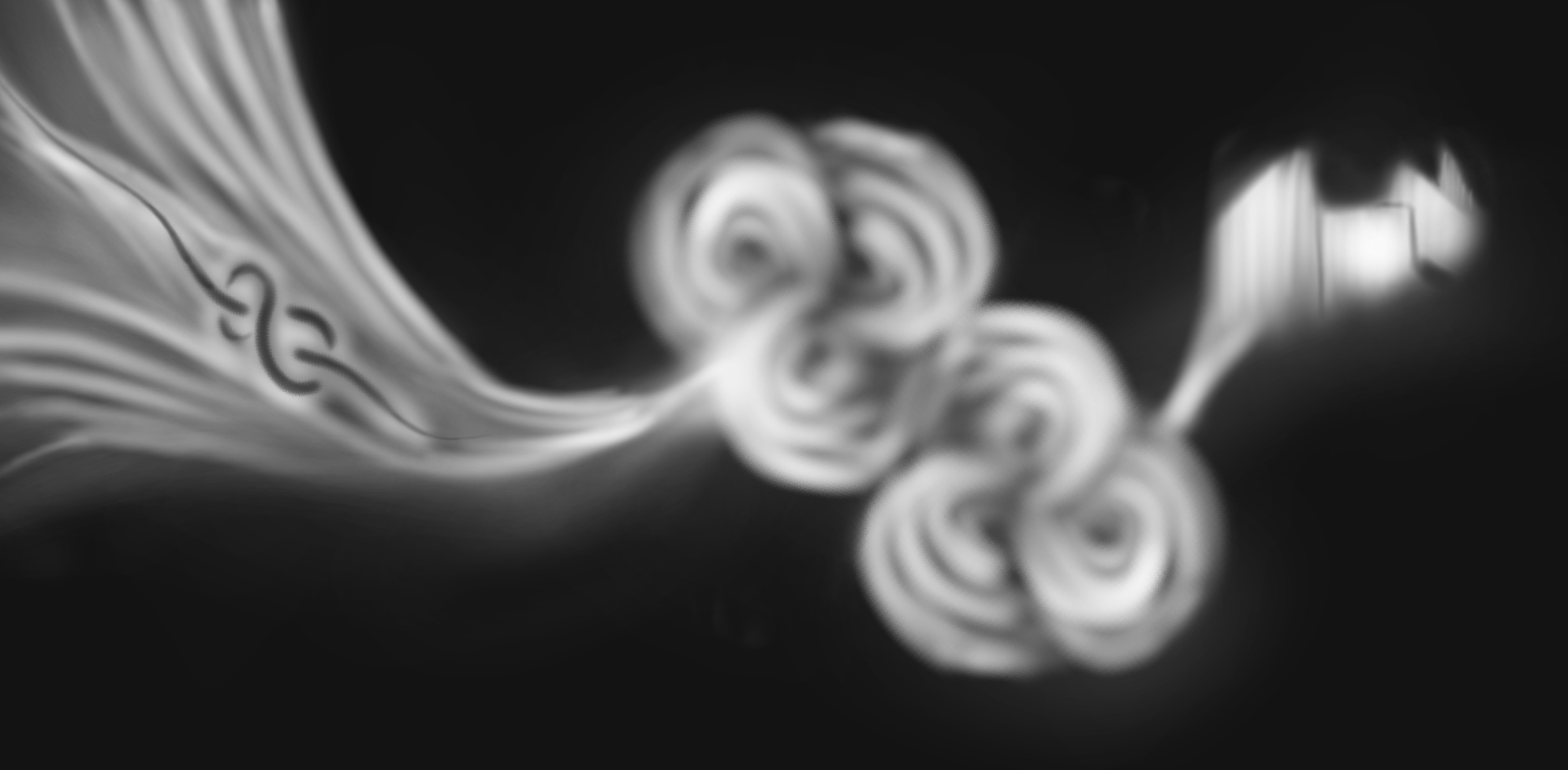}
    \caption{A space formed from three different geometries\ldots}
    \label{fig:Abstract three manifold}
    \end{figure}
    
    \begin{figure}[H]
       \centering
    \includegraphics[width=.9\linewidth]{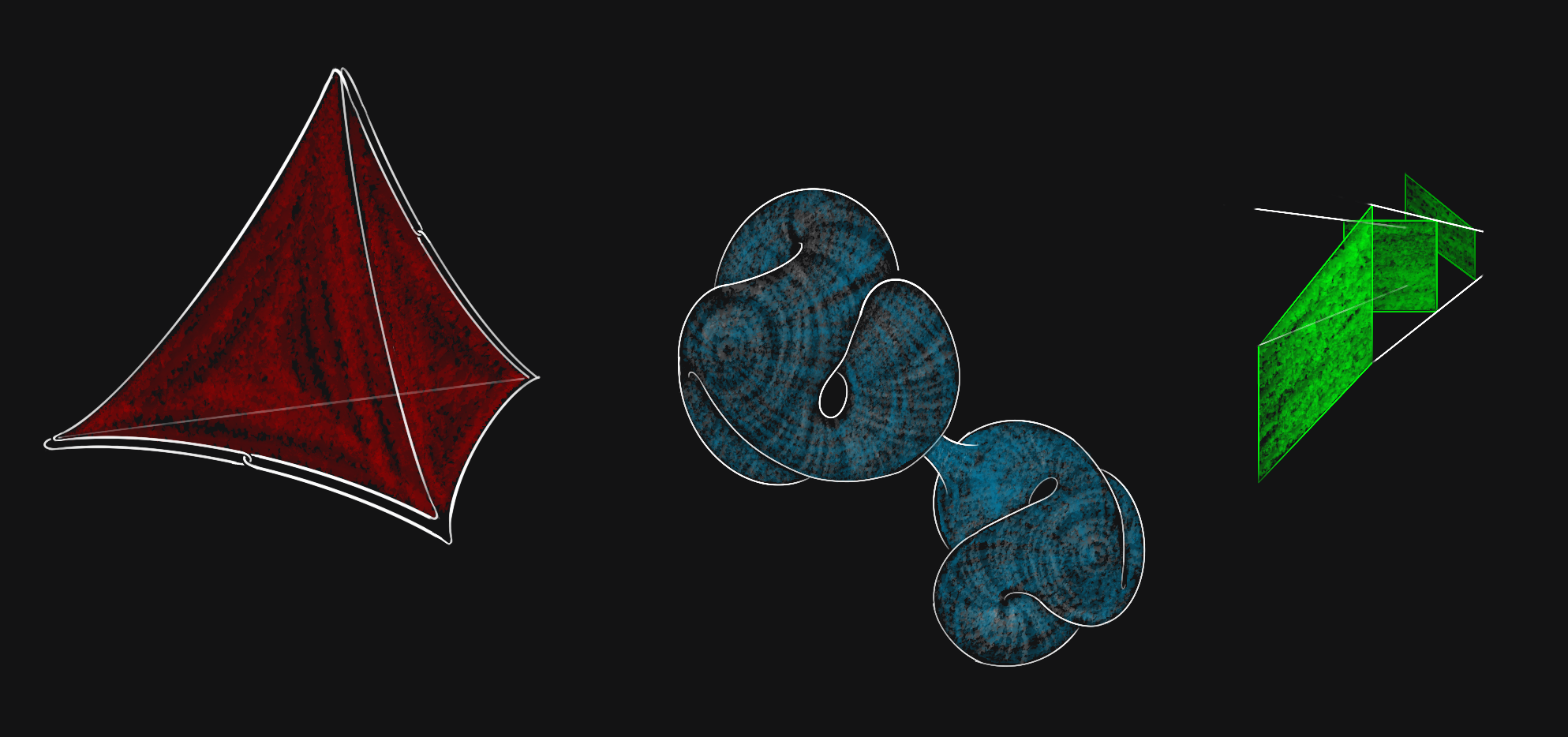}
    \caption{And its geometric decomposition\protect\footnotemark}
    \label{fig:Geometrizationdecomposition}
\end{figure}

\footnotetext{This figure incorporates some artistic license as the complement of a figure-eight knot is non-compact, so the space as depicted does not satisfy the hypotheses of the Geometrization conjecture. However, it seemed necessary to include at least one depiction of a hyperbolic knot complement while discussing Thurston's work. The image of the knot is adapted from the lower figure on page 5 of \cite{thurston1979geometry}. The central component of this manifold is a $\mathbb{S}^2 \times \mathbb{R}$ structure on the connected sum of two copies of $\mathbb{RP}^3$. Since there is not a good way to draw three-dimensional projective space, what is shown is a connected sum of two Boy's surfaces, which are immersions of $\mathbb{RP}^2$ in $\mathbb{R}^3$. The right-most piece has Nil geometry, which is depicted by a portion of the Cayley graph of the Heisenberg group.}

\section*{The Ricci flow}

Around this time, Richard Hamilton proposed an ambitious program to attack both of the Geometrization and the Poincar\'e conjectures \cite{hamilton1982three}. He had been studying work by James Eells and Joseph Sampson~\cite{eells1964harmonic}, which used ideas about heat flows to find harmonic maps and thought it might be possible to use a similar approach for these problems.\footnote{Strictly speaking, the Ricci flow is not a (non-linear) heat flow. In particular, it is diffeomorphism-invariant, which induces zero terms in its symbol. As a result of this, establishing the existence for the flow was a major accomplishment. Hamilton did so via a technical argument ~\cite{hamilton1982three} but a few years later Dennis DeTurck found a much simpler argument by conjugating the Ricci flow to obtain a strictly parabolic flow \cite{deturck1981existence}.}
  He defined an evolution equation, known as the Ricci flow, which would deform the shape of a space and hopefully allow its curvature to dissipate throughout the space. If the space was simply connected, he hoped the geometry would evolve to that of a round sphere, which would establish the result. 
 
Hamilton and others were able to make partial progress towards this goal. In particular, Hamilton showed that starting with a simply connected three-dimensional space whose Ricci curvature is positive, it would indeed evolve to a round sphere. This was a weaker version of the Poincar\'e conjecture, but a major proof a concept for the power of the Ricci flow. Furthermore, he was nearly able\footnote{There were some missing steps dealing with metrics of mixed curvature on genus zero surfaces which were established by Bennett Chow \cite{chow1991ricci} and Xiuxiong Chen, Peng Lu, and Gang Tian \cite{chen2006note}.} to find a new proof of the uniformization theorem using Ricci flow. Although the Ricci flow was a promising approach to the Poincar\'e conjecture, the flow would encounter \emph{singularities}, which are times where the space would either collapse to a point\footnote{The Ricci flow can also collapse in other ways. For instance, if we consider the Cartesian product of a round sphere with a unit circle (with the product metric), Ricci flow will collapse the space $\mathbb{S}^2\times \mathbb{S}^1$ to the unit circle.} or violently rip itself apart. Hamilton understood the cases when the space would shrink to a point, but was was not able to control the geometry in the latter case, and this presented a fundamental obstacle to finishing the proof.  

\begin{figure}[H]
    \centering
\includegraphics[width=.9\linewidth]{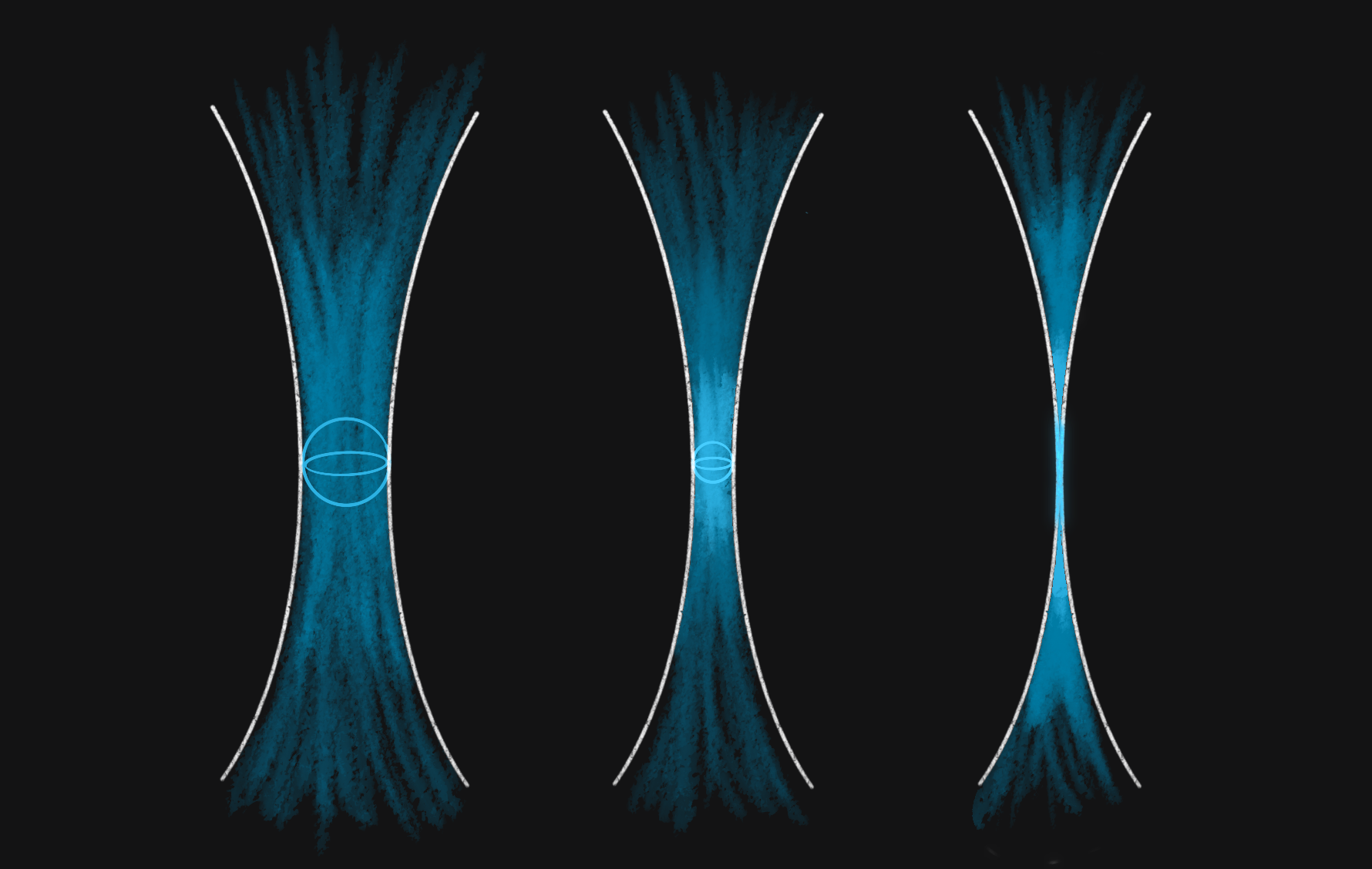}
    \caption{The Ricci flow approaching a singularity}
    \label{fig:A singularity forming}
\end{figure}

\section*{Perelman's breakthrough}

In 2000, the Clay Millenium institute listed seven major open problems in mathematics and provided a \$1,000,000 prize for a solution to any one of them. The Poincar{\'e} conjecture was chosen as one of the problems and despite the best efforts of many mathematicians to solve it, no one expected to see a solution in the near future.

However, just two years later, Grigori Perelman posted a preprint to the arXiv \cite{perelman2002entropy} with several major breakthroughs in the Ricci flow. In particular, he showed that by coupling the Ricci flow with a heat equation running in reverse, there is a non-decreasing quantity (which corresponds to the \emph{Fisher information} of the heat distribution\footnote{For background on the Fisher information, we refer the reader to Chapter 20 of Cedric Villani's textbook Optimal Transport, old and new \cite{villani2009optimal}.}) which can be used to control the geometry of its solutions. Perelman claimed that this idea could be used to solve both the Poincar\'e and Geometrization conjectures, and gave a short sketch of how to do this. This paper took the geometry community by storm, and the experts began to study it to understand what Perelman had found.

\begin{figure}[H]
    \centering
\includegraphics[width=.9\linewidth]{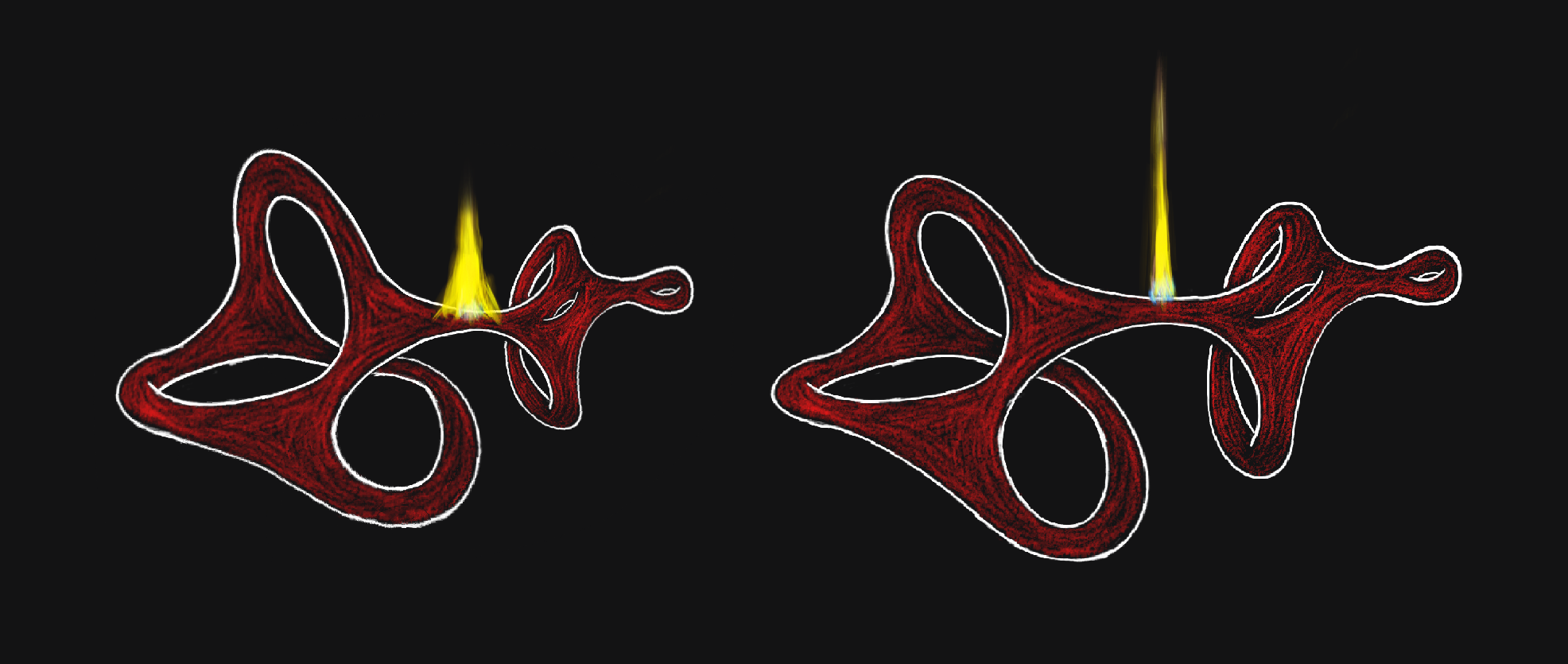}
    \caption{An evolving surface with a backwards heat flow}
    \label{fig:Backwardsheatequation}
\end{figure}

Over the eight months, two more preprints \cite{perelman2003ricci, perelman2003finite} followed the first and provided more details. In particular, Perelman combined the ideas from the first paper with a geometric process known as surgery to excise regions of space when they started to tear apart and replace them with pieces with better behavior.\footnote{Ricci flow with surgery was invented by Hamilton in 1993 \cite{hamilton1993formations,hamilton1997four}, but he was unable to use it to handle the singularities in the three dimensional case.} By combining the Ricci flow with surgery, Perelman was able to show that any simply-connected three-dimensional space would converge to a round sphere (or possibly a connected sum of several round spheres), and thus is topologically equivalent to the standard sphere. For more general three-dimensional spaces, he proved that after a large amount of time, Ricci flow with surgery would decompose the space into several pieces whose geometric structure was well understood.\footnote{Ricci flow with surgery does not necessarily converge to one of Thurston's geometries on each connected component. Indeed, several of the geometries collapse under the flow, so are not even fixed points. Instead, what occurs is that if we perform the surgeries in an effective way, there are only finitely many surgeries \cite{bamler2018long}, and after the last one the space decomposes into finitely many pieces, all of which were previously known to satisfy the Geometrization conjecture \cite{calegari2020ricci}.} This established both the Poincar{\'e} conjecture as well as the Geometrization Conjecture.

\begin{figure}[H]
    \centering
\includegraphics[width=.9\linewidth]{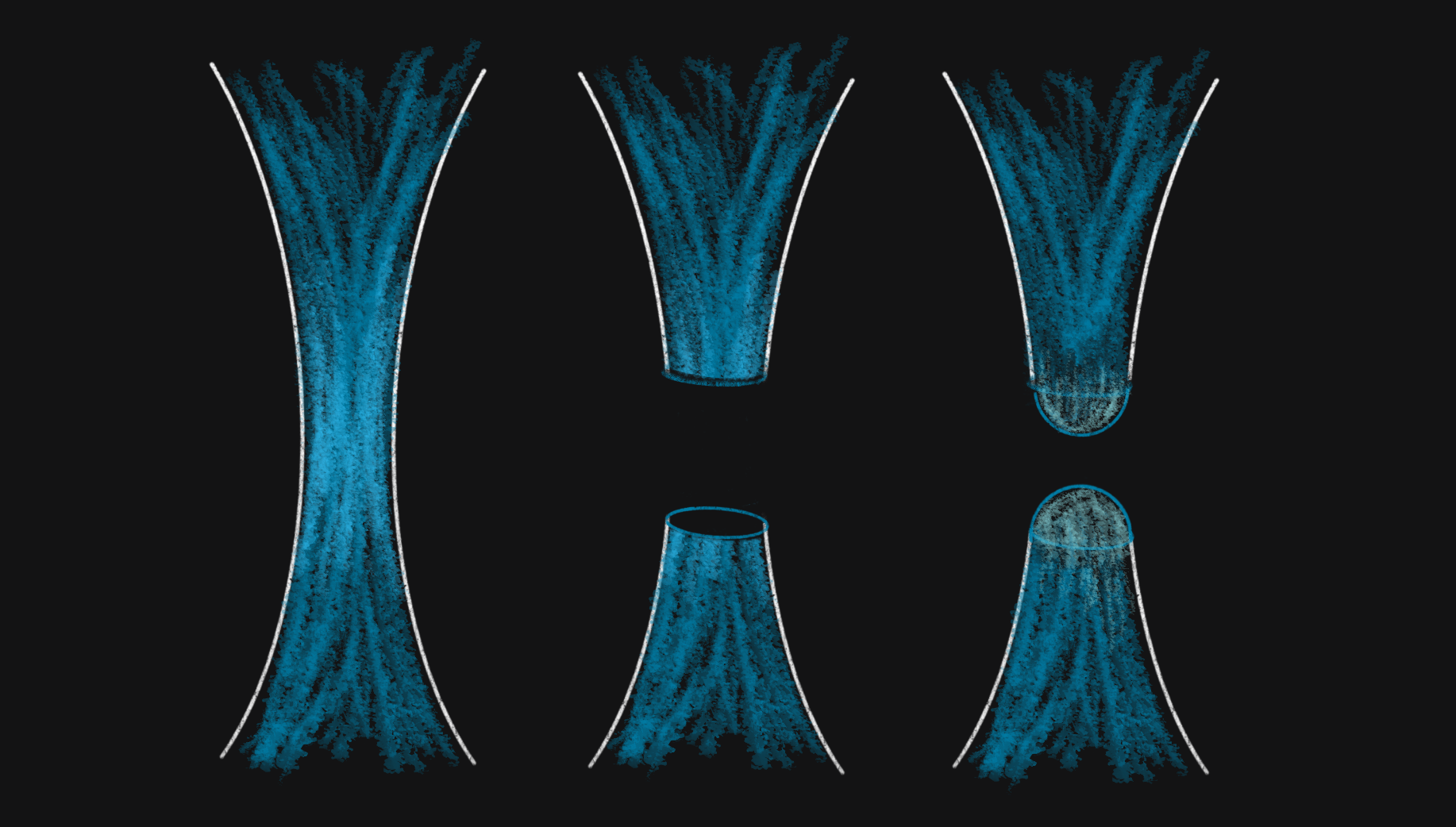}
    \caption{Using surgery to tame singularities}
    \label{fig:Singularities}
\end{figure}

Perelman had worked for nearly a decade in relative isolation, and his proof was the crowning achievement in a long line of work stretching back hundreds of years. The result sparked an enormous amount of interest in geometric analysis and the Ricci flow in particular. However, it took several years (and an authorship scandal) for the mathematical community to accept the work as correct. In 2006, Perelman was awarded a Fields medal for his work, but declined the award and later turned down the Millenium Prize, as well. He has since withdrawn from mathematics.

\chapter*{What is the Ricci flow?}
\addcontentsline{toc}{chapter}{What is the Ricci flow?}

The goal of this chapter is to provide a working definition for the Ricci flow with some intuition for its behavior. Heuristically, the Ricci flow is a ``heat flow of curvature" and we will try to explain what this means. To do so, we have divided this chapter into four sections.
\begin{itemize}
    \item[] \hyperref[Heat equation section]{The heat equation and its behavior} \hfill 14
    \item[]  \hyperref[Curvature section]{The curvature of space} \hfill 21
    \item[] \hyperref[Ricci Flow section]{The Ricci flow as a geometric heat flow} \hfill 28
    \item[] \hyperref[Distinctions section]{Differences between Ricci flow and the standard heat equation} \hfill 31
\end{itemize}

\section*{The Heat Equation}
\label{Heat equation section}

The heat equation models how the temperature of a region changes in time. More precisely, we consider a function $u(\mathbf{x},t)$ which returns the temperature of a point\footnote{Here we are using boldface $\mathbf{x}$ to distinguish from the $x$-coordinate.} $\mathbf{x} \in \mathbb{R}^n$ at a time $t$.
Ignoring what happens at the boundary,\footnote{In the context of the Ricci flow, we will be dealing with compact manifolds without boundary so this is not an issue.} we say that the function $u(\mathbf{x},t)$ solves the heat equation if it satisfies 
\begin{equation} \label{Heat equation}
    \frac{ \partial u(\mathbf{x},t)}{\partial t} = \Delta u(\mathbf{x},t).
\end{equation} 
This equation was introduced by Joseph Fourier in 1822, and played a foundational role in the development of Fourier analysis. We will not provide a physical derivation for why heat distributions tend to obey this equation, but instead try to understand the behavior of its solutions.

\begin{figure}
    \centering
    \includegraphics[width=.9\linewidth]{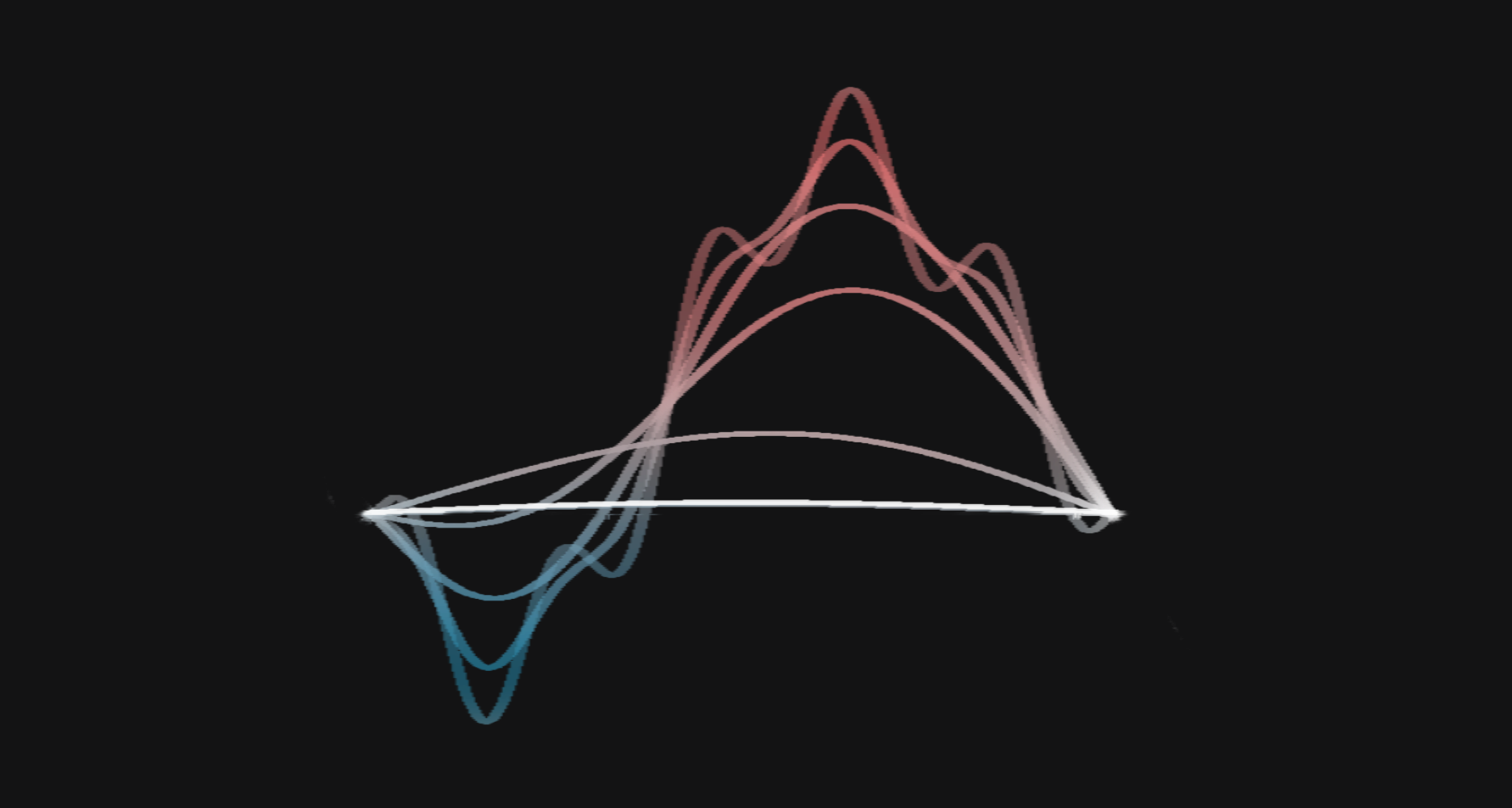}
    \caption{A heat distribution flowing to equilibrium\protect\footnotemark}
    \label{fig:heatflow}
\end{figure}

\footnotetext{In this diagram, we have imposed Dirichlet boundary conditions on the solution. The time slices were chosen to show qualitative behavior of the solution, and were not taken uniformly.}

\subsection*{The Laplacian}
Before trying to understand the behavior of solutions to Equation \ref{Heat equation}, let us first try to understand each of its terms. The left hand side of Equation \ref{Heat equation} is the partial derivative of $u$ with respect to $t$, which describes how the temperature at a point $\mathbf{x}$ evolves as time goes by. Therefore, we must understand the term $\Delta u$, which is the Laplacian of $u$. In vector calculus, this is often defined in the following way:  \[ \Delta u =
\textrm{div}(\textrm{\,grad } u ).\] However, there are several other perspectives on the Laplacian that will tie in more naturally to Ricci flow. One useful interpretation of the Laplacian is to consider the (three-dimensional) Hessian matrix
\[  \Hess(u) =
\begin{bmatrix}
    \frac{\partial^2 }{\partial x^2} u & \frac{ \partial }{\partial x}\frac{ \partial }{\partial y} u & \frac{ \partial }{\partial x}\frac{ \partial }{\partial z} u  \\
     \frac{ \partial }{\partial y}\frac{ \partial }{\partial x} u & \frac{\partial^2 }{\partial y^2}  u& \frac{ \partial }{\partial y}\frac{ \partial }{\partial z} u   \\
    \frac{ \partial }{\partial z}\frac{ \partial }{\partial x}  u& \frac{ \partial }{\partial z}\frac{ \partial }{\partial y} u & \frac{\partial^2 }{\partial z^2} u
\end{bmatrix}
\]
and note that the Laplacian is the trace of this matrix: \[\Delta u = \frac{\partial^2 }{\partial x^2} u + \frac{\partial^2 }{\partial y^2} u +\frac{\partial^2 }{\partial z^2} u. \]
Later on, we will see that the Ricci curvature is the trace of the Riemann curvature tensor, so this identity makes the connection between the Ricci curvature and the Laplacian fairly natural.

\begin{figure}[H]
    \centering
\includegraphics[width=.9\linewidth]{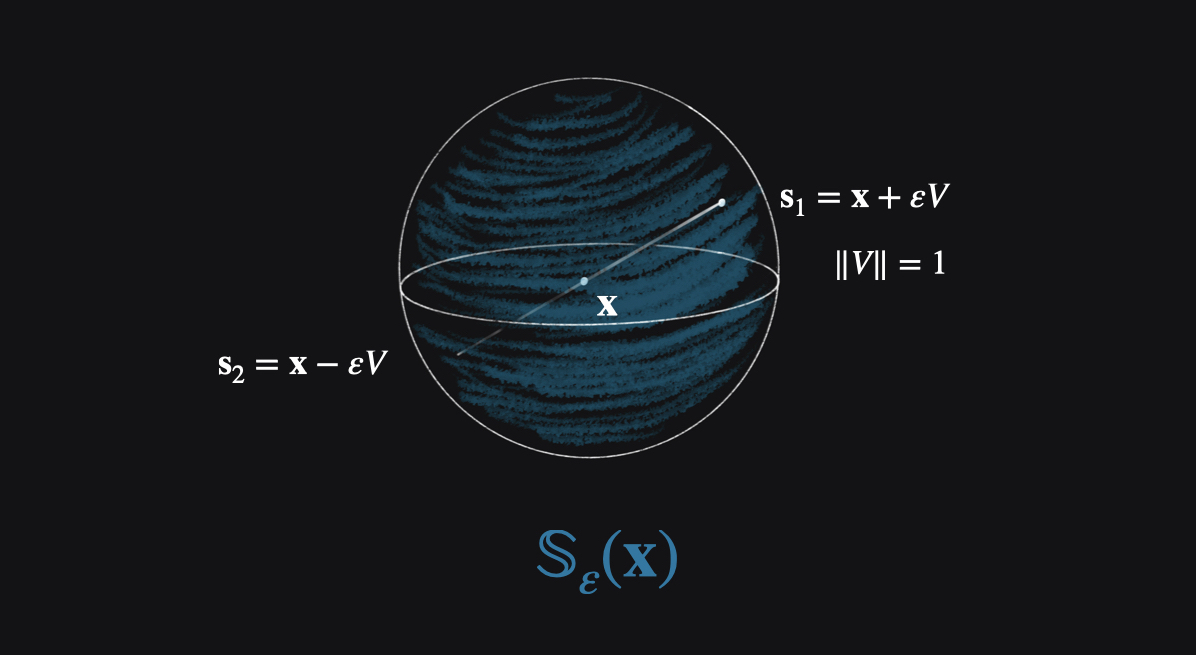}

    \caption{The Laplacian as an integral}
    \label{fig:Laplacian as integral}
\end{figure}
		\begin{center}

\end{center}

In order to define the Ricci flow using the minimal amount of Riemannian geometry, there is a third perspective on the Laplacian that is even more useful. We consider a point $\mathbf{x}$ and a small sphere of radius $\varepsilon$ around $\mathbf{x}$, which we denote $\mathbb{S}_\varepsilon(\mathbf{x})$. 
One can show that the Laplacian is proportional to the difference between $f(\mathbf{x})$ and the average value of $f$ on $\mathbb{S}_\varepsilon(\mathbf{x})$:

\begin{equation} \label{Laplace as Average}
    \Delta f(\mathbf{x})=\lim _{\varepsilon \rightarrow 0+} \frac{2 n}{\varepsilon^{2}} \frac{1}{\omega\left(\mathbb{S}_{\varepsilon}\right)} \int_{\mathbb{S}_{\varepsilon}(\mathbf{x})}(f(\mathbf{s})-f(\mathbf{x})) \, \mathrm{d} \omega(\mathbf{s}).
\end{equation} 

Here, $\omega\left(\mathbb{S}_{\varepsilon}\right)$ is the surface area of a sphere of radius $\varepsilon$ in Euclidean $n$-space, which is given explicitly by  \[ \omega(\mathbb{S}_\varepsilon)=\frac{2\pi^\frac{n}{2}}{\Gamma\left(\frac{n}{2} \right)} \varepsilon^{n-1}. \]

From this we can see that $\Delta u$ is the average of second derivatives, which is a perspective that will help us understand Ricci curvature. It is worth convincing yourself that the right-hand side of \eqref{Laplace as Average} really does have something to do with second derivatives.\footnote{ Figure \ref{fig:Laplacian as integral} is meant to serve as a hint for why this is the case. As an additional hint, given a unit vector $V$, what is \[ \lim_{\varepsilon \to 0}\frac{ f(\mathbf{x}+\varepsilon V) - f(\mathbf{x}) + f(\mathbf{x}-\varepsilon V)-f(\mathbf{x})}{\varepsilon^2}? \]}

\subsection*{Convergence to equilibrium}
%The heat equation is a partial differential equation given by the expression
%$$\frac{ \partial u}{\partial t} = \Delta u.$$

One hallmark of the heat equation is that its solutions tend to converge to an equilibrium. That is to say, given a function $u(\mathbf{x},t)$ which solves the heat equation, we expect that \[ \lim_{t \to \infty} u(\mathbf{x},t) = C \]
where $C$ is a constant which is independent of the point.
Intuitively, if we have a warm cup of coffee outside in a cold day, this states that the heat from the coffee will dissipate into the air and the temperature will converge to that of the outside environment.

There are various ways to make this idea into a precise mathematical statement, and we will mention two that are relevant for understanding the Ricci flow.

\subsubsection*{The maximum principle}

One of the most fundamental tools that one has for studying heat equations is the \emph{maximum principle.} To give a simple illustration of this principle, we will show that the hot spots of $u(\mathbf{x},t)$ tend to cool down while the cold spots tend to warm up. To do so, we suppose the temperature $u$ is maximized at the point $\mathbf{x}_0$ (in the interior of the domain) at some time $t$. Then, since $u(\mathbf{x}_0,t)$ is maximized, the second derivative test implies that the Hessian of $u$ at $\mathbf{x}_0$ must be non-positive definite (i.e., have all non-positive eigenvalues). This implies that $ \Delta u(\mathbf{x}_0,t)$ must also be non-positive. The heat equation then tells us that \[ \frac{ \partial u(\mathbf{x}_0,t)}{\partial t} \leq 0, \]
which is to say that the temperature at $\mathbf{x}_0$ is decreasing.
If we consider the coldest points, the same argument shows that they warm up under the heat flow. 

 This type of argument, where one considers a point which maximizes some function and applies the second derivative test (or a more sophisticated version of it) is endlessly useful, and plays an essential role in the analysis of heat-type equations. For example, one can use this argument to prove the Li-Yau estimate \cite{li1986parabolic}, which states that positive solutions to the heat equation on a space with non-negative Ricci curvature satisfies the inequality
\begin{equation*}
    \Delta \ln u(\mathbf{x},t) \geq  - \frac{n}{t}.
\end{equation*}
This particular estimate plays a fundamental role in geometric analysis, including in the study of the Ricci flow.

%In order to show that the temperature of the hottest point is strictly decreasing, we would need a stronger version of the maximum principle. The maximum principle alone also doesn't show that the heat converges to a constant as $t$ goes to infinity, but for reasonable initial conditions, this is indeed the case.\footnote{For instance, using Fourier analysis it is possible to show that the heat will converge to a constant if the initial data is $W^2([0,1])$ (or even less regular) with periodic boundary conditions.}

\subsubsection*{Entropy and energy}

Another fundamental tool for heat-type equations is to find quantities like energy or entropy which either increase or decrease along the flow. By doing so, one can better understand the behavior of the heat equation.

\begin{figure}[H]
    \centering
  \includegraphics[width=.9\linewidth]{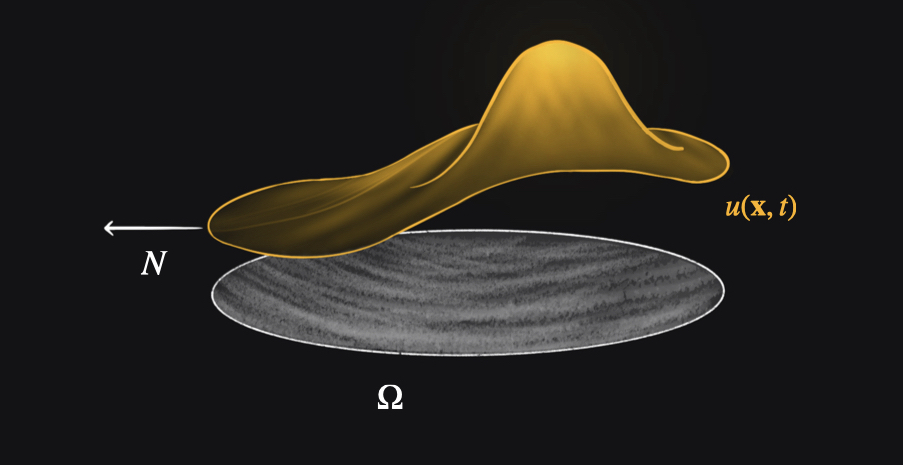}
    \caption{A perfectly insulated domain}
    \label{fig:Neumann boundary conditions}
\end{figure}

To show this idea in practice, we consider the heat equation on a smooth bounded domain $\Omega$ and assume that the domain is perfectly insulated from its surroundings. %so that there is no distant heat source warming the domain. 
From a mathematical perspective, this assumption is equivalent to requiring that the temperature at the boundary of the domain satisfies \[ \nabla_N u(\mathbf{x},t) =0\]
where $N$ is the outer normal vector, which is known as \emph{Neumann boundary conditions}.\footnote{For those who have taken a class on PDEs, it is a good exercise to show that the total heat is conserved under Neumann boundary conditions.}

We then consider the quantity \begin{equation} \label{L2 energy}
 E(t) =  \int_\Omega  u(\mathbf{x},t)^2 \, d \mathbf{x},
\end{equation}
which is often called the \emph{energy}. However, the name should not be taken too literally because it doesn't have a direct physical interpretation. Intuitively, it measures how concentrated the heat distribution is. In other words, if $u(\mathbf{x},t)$ is very large in some region $\mathbf{x}$, $E(t)$ will also be large. On the other hand, if the distribution is very diffuse, $E(t)$ will be much smaller. We can compute how $E(t)$ evolves in time as follows.

\begin{eqnarray}
\frac{d E}{dt} & = &  \int_\Omega  \frac{\partial }{ \partial t} u(\mathbf{x},t)^2 \, d \mathbf{x} \\
% & = & \int_\Omega  2 \frac{\partial u(\mathbf{x},t) }{ \partial t}  \cdot u(\mathbf{x},t) \, d \mathbf{x} \\
& = & \int_\Omega  2 \left ( \Delta u(\mathbf{x},t) \right) \cdot u(\mathbf{x},t) \, d \mathbf{x}.
\end{eqnarray}

Using integration by parts and the generalized Stokes' theorem, it is possible to simplify this expression.

\begin{eqnarray}
\frac{d E}{dt} & = &  - \int_\Omega 2  \| \nabla u(\mathbf{x},t) \|^2  \, d \mathbf{x} +  \int_{\partial \Omega } 2 u(\mathbf{x},t) \left \langle \nabla u(\mathbf{x},t), N(\mathbf{x}) \right \rangle \, d \mathbf{x}.
%& = &  - \int_\Omega 2  \| \nabla u(\mathbf{x},t) \|^2  \, d \mathbf{x},
\end{eqnarray}

However, the boundary conditions implies that the second term vanishes.

\[ \frac{d E}{dt} = - \int_\Omega 2  \| \nabla u(\mathbf{x},t) \|^2  \, d \mathbf{x} < 0. \]

As such, the energy decreases along the flow, which shows that the heat distribution tends to spread out through space. With a bit more effort, it is also possible to show that this quantity is convex.

Another quantity that we can consider is \[ S = - \int_\Omega u(\mathbf{x},t) \ln \left(u(\mathbf{x},t) \right ) \, d \mathbf{x}, \]
which is better known as the \emph{entropy}.\footnote{Here, we need to assume that the temperature is everywhere positive for this quantity to be well-defined.} The entropy is a measure of how disordered a configuration is. More precisely, it is the amount of information (measured in \emph{nats}) that the heat distribution contains relative to the equilibrium state. However, entropy is a notoriously tricky concept to conceptualize, so for our purposes we can simply define it using the integral above. We can use the same idea from before to calculate the evolution of $S$. Doing so, we find the following:

\begin{eqnarray}\label{Entropy and information}
   \frac{d S}{d t} & = &  \int_\Omega \frac{\|\nabla u\|^2}{u}  \, d \mathbf{x} > 0
\end{eqnarray}

\begin{comment}
\begin{eqnarray}
   \frac{d S}{d t} &= & - \int_\Omega \frac{\partial}{\partial t} (u \log u) d \mathbf{x} \\
   & = & - \int_\Omega (\Delta u) \log u +\frac{u}{u} \Delta u \, d \mathbf{x} \\
   & = &  \int_\Omega \frac{\|\nabla u\|^2}{u}  \, d \mathbf{x} > 0
\end{eqnarray}
  \end{comment}
  
    The quantity on the right hand side of this equation is known as the Fisher information, and is positive. As such,
Equation \ref{Entropy and information} is a mathematical version of the second law of thermodynamics, which states that a closed system tends to go from an orderly configuration to a disordered state. These types of quantities play a central role in the Ricci flow and one of Perelman's most important breakthroughs was to find a version of ``entropy" for the Ricci flow.\footnote{In fact, his first paper in the trilogy is titled ``The entropy formula for the Ricci flow and its geometric applications" \cite{perelman2002entropy}. However, the precise quantity is more similar to a Fisher information rather than an entropy.}

\section*{Curvature}
\label{Curvature section}
In order to discuss the Ricci flow, it is first necessary to discuss the notion of curvature. Unfortunately, rigorously defining curvature requires some background knowledge in Riemannian geometry and a more in-depth discussion of extrinsic versus intrinsic geometry, which are both outside the scope of this introduction. As such, in this section we will provide an intuitive notion of curvature, without trying to be overly precise.

\subsection*{Sectional curvature}

It is not immediately obvious what curvature is, especially from an intrinsic perspective where the space of interest does not live in some ambient Euclidean space. However, one useful bit of intuition is that Euclidean space is flat, without any bumps or valleys. As such, one way to formalize the notion of curvature is to compare the geometry of a given space with that of Euclidean space. To do so, the simplest approach is to consider triangles.% in Euclidean space to triangles in curved space.

\begin{figure}[H]
    \centering
  \includegraphics[width=.9\linewidth]{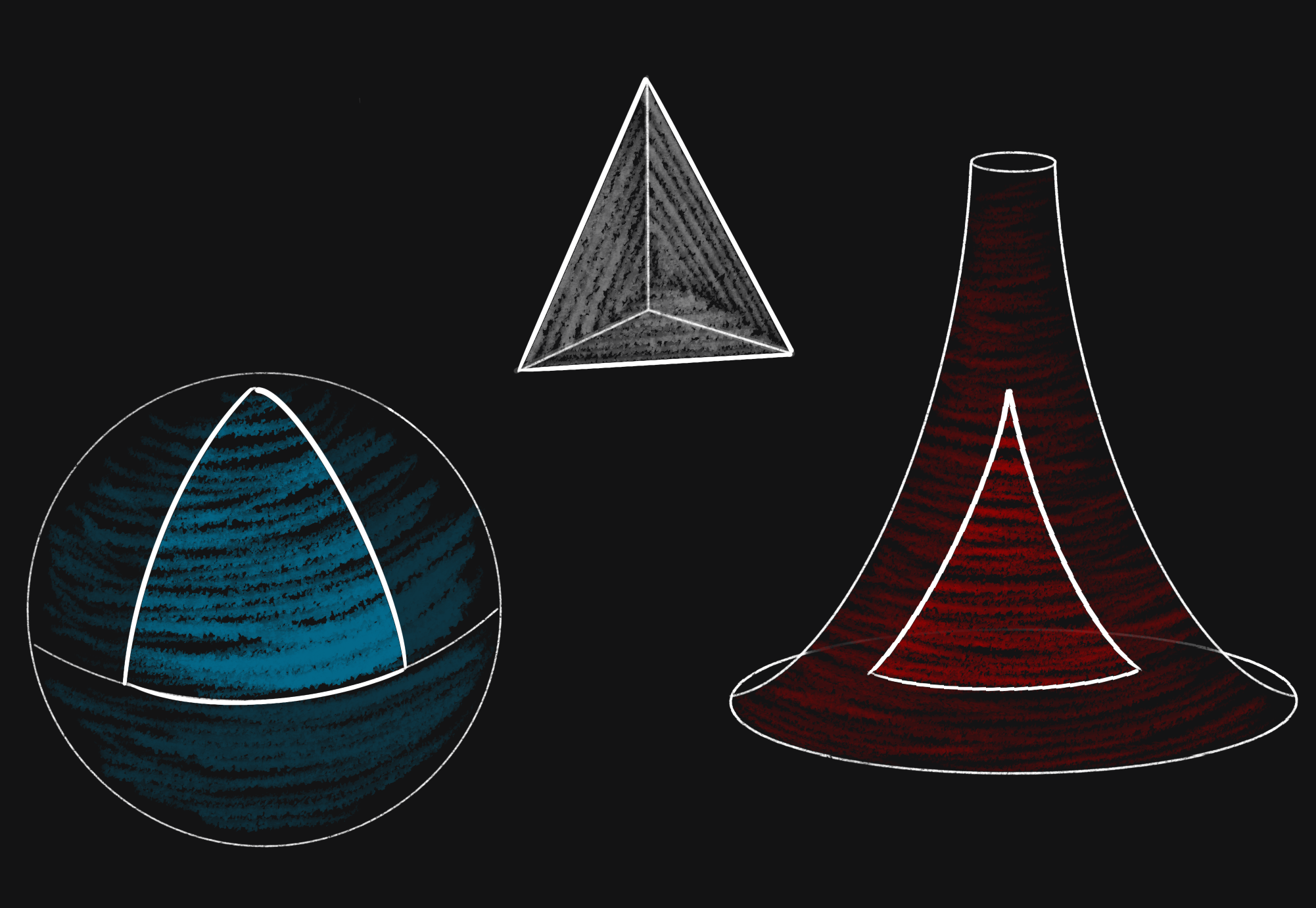}
    \caption{Triangles in spherical, flat and hyperbolic space}
    \label{fig:Triangle comparison}
\end{figure}

As can be seen in Figure \ref{fig:Triangle comparison}, triangles on a sphere\footnote{More precisely, we are considering geodesic triangles, where each side is length-minimizing.} appear ``fatter" than triangles in flat space. There are a few ways to make this precise. For instance, we can consider the sum of the angles of a spherical triangle. It turns out the sum is greater than $\pi$, and actually depends on the area of the triangle.\footnote{For readers who have studied spherical trigonometry, this will be a familiar fact. On a two dimensional sphere, it is actually a consequence of the Gauss-Bonnet theorem and for those who are familiar with the differential geometry of curves and surfaces, it is a good exercise to compute the exact formula for the sum of the angles.}
 Apart from the angles, we can also see from the picture that the sides of a spherical triangle seem to bow away from the other sides. This bowing occurs in any space of positive curvature, and is one of the characteristic features of positive sectional curvature. To understand this, it helps to picture the sides as turning towards each other, which is why they appear to bow outward. Here, the meaning of the word ``turning" is  somewhat imprecise, but hopefully the picture makes it clear what we mean.

Before we use this idea to discuss curvature, it is worthwhile to note that round spheres are quite special, in that every point looks like every other point.\footnote{Round spheres are examples of  \emph{symmetric spaces}.} However, the earth is not exactly a sphere; the geometry of Mount Everest looks very different from that of the Great Plains. This will be true for most of the spaces we are interested in, so we want some way to define curvature on spaces which are not so symmetric. To do so, instead of using a giant triangle as we did on the sphere, we use very small ones.

To do this, we consider a point in our space and two tangent vectors $X$ and $Y$. Then we consider a triangle which has one vertex at $p$ and which has two sides of length $\varepsilon$ in the directions of $X$ and $Y$. We will not be too precise about what it means to make a triangle where the sides ``head in a direction," because it takes some work to define.\footnote{More precisely, we take a tangent vector $X$ and consider the geodesic $\gamma(t) = \exp_p(Xt)$.} On the sphere, we saw that the triangle will turned in on itself as a result of the positive curvature. With this in mind, we consider the length of the third edge of the triangle, which we denote $L(\varepsilon)$.

\begin{figure}[H]
    \centering
\includegraphics[width=.9\linewidth]{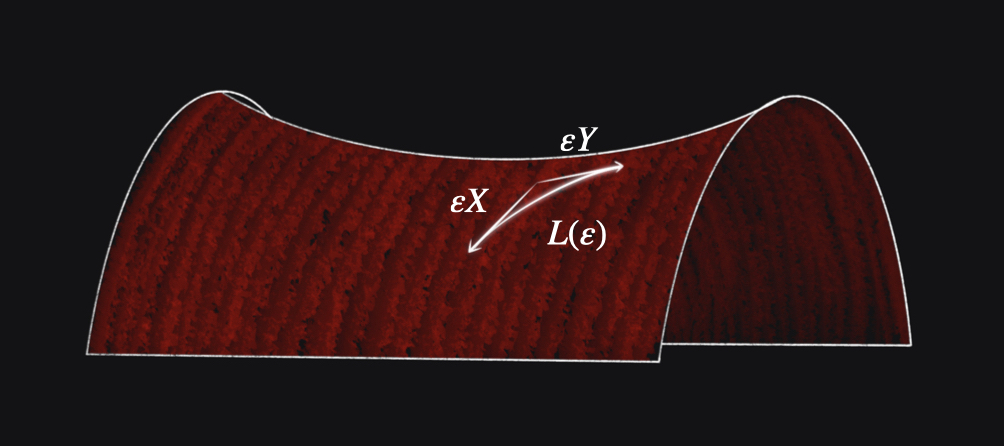}
    \caption{A heuristic diagram for $L(\varepsilon)$}
    \label{fig:L epsilon diagram}
\end{figure}

The curvature will be positive whenever $L(\varepsilon)$ is smaller than it would be for a corresponding triangle\footnote{Here, a \emph{corresponding triangle} means a triangle $\triangle PQR$ whose angle $\angle RPQ$ is the equal to the angle between $X$ and $Y$ and where the sides $\overline{PQ}$ and $\overline{PR}$ have length $\varepsilon$.} in Euclidean space and negative if $L(\varepsilon)$ is greater than a triangle in Euclidean space. To make this idea more formal, we consider the Taylor polynomial of $L(\varepsilon)$ in terms of $\varepsilon$. Doing so, we find the following:
\begin{equation} \label{Taylor expansion of L}
    L(\varepsilon)=\varepsilon\|X-Y\|\left(1-\frac{1}{12} K(X, Y)(1+\langle X, Y\rangle) \varepsilon^{2}\right)+O\left(\varepsilon^{4}\right)
\end{equation}

In this expression, $K(X,Y)$ is defined to be the sectional curvature of the tangent plane spanned by $X$ and $Y$. Note that in flat space, $K(X,Y) \equiv 0$ so this matches with our intuition. On the other hand, the unit sphere has constant positive sectional curvature (i.e., $K(X,Y) \equiv 1$). Since $K(X,Y)$ is positive, the third side of the triangle will be shorter than the corresponding side in flat space. Conversely, whenever the sectional curvature is negative, the third side of the triangle will be longer than that of the corresponding Euclidean triangle.

\subsection*{Ricci curvature}

\label{Ricci curvature section}

While the sectional curvature contains all of the curvature information about a space, it is a very complicated object. The Ricci curvature is a coarser invariant than the sectional curvature,\footnote{In dimensions two and three, the Ricci curvature determines the sectional curvature completely, but this is not the case in higher dimensions.} but conveys important geometric information that is essential in many applications. To define it, we consider a unit vector $X$ and define the Ricci curvature\footnote{The standard definition of the Ricci curvature is as the contraction of the Riemann curvature tensor along second and last indices. In other words, given an orthonormal frame $\{e_i\}_{i=1}^n$ and two vector fields $X$ and $Z$,  \[\Ric(X,Z) = \sum_{i=1}^n \langle R(X,e_i)Z, e_i \rangle \] Here, $R(\cdot, \cdot) \cdot$  is the Riemannian curvature tensor, which is defined as  \[R(X,Y)Z = \nabla_Y \nabla_X Z - \nabla_X \nabla_Y Z - \nabla_{[X,Y]} Z. \] where $\nabla_X$ denotes covariant derivative in the $X$ direction with respect to the Levi-Civita connection. We used Equation \eqref{Ricci as average} to avoid having to define the notions of Riemannian curvature and covariant derivatives.} $\Ric(X,X)$ to be $(n-1)$ times the average of all of the sectional curvatures of tangent planes containing $X$. In other words, the Ricci curvature satisfies the identity\footnote{Note that the pairs $(X,Y)$ and $(X,-Y)$ span the same plane, which is why there is a factor of $\frac{1}{2}$ in Equation \ref{Ricci as average}.} 
\begin{equation} \label{Ricci as average}
     \Ric(X,X) = \frac{1}{2} \frac{(n-1)}{\omega\left(\mathbb{S}^{n-2}\right)} \oint_{\|Y \| = 1 \textrm{ and } X \perp Y} K(X,Y) \, \mathrm{d} \mathbb{S}^{n-2}(Y),
\end{equation}
where $\omega \left(\mathbb{S}^{n-2}\right)$ is the surface area of the $(n-2)$-dimensional sphere. Initially, it might seem a bit strange that the vector $X$ appears twice as an argument for the Ricci curvature, but we will just treat this as a convention without going into detail about why this is the case.\footnote{At each point, the Ricci curvature is a symmetric bilinear form, which is why there are two copies of $X$. To obtain the bilinear form from the average of sectional curvatures, one uses the polarization identity \[ \Ric(X,Y) = \frac{1}{2} \left(\Ric(X+Y,X+Y) - \Ric(X,X) - \Ric(Y,Y)\right) \].}

Equation \eqref{Ricci as average} gives a concise definition for the Ricci curvature, but does not provide any intuition for what the Ricci curvature actually is. To do so, it is possible to draw a picture for the Ricci curvature which is similar to the one for sectional curvature but uses cones rather than triangles.\footnote{The idea of Ricci curvature as the distortion of narrow geodesic cones is taken from Chapter 14 of Villani's textbook on optimal transport \cite{villani2009optimal}.}

We consider a point $p$, a unit vector $X$ at $p$ and take a short segment of length $\varepsilon$ in the $X$ direction.\footnote{More precisely, we consider a geodesic $\gamma$ of length $\varepsilon$ with $\dot \gamma(0) = X$.} Finally, we put a narrow circular cone with vertex $p$ around the segment (here, the cone being narrow means that the cone angle $\theta$ satisfies $\tan(\theta) = \varepsilon$). We then consider the area of the base of the cone, which we denote $A(\varepsilon)$.

\begin{figure}[H]
    \centering
\includegraphics[width=.9\linewidth]{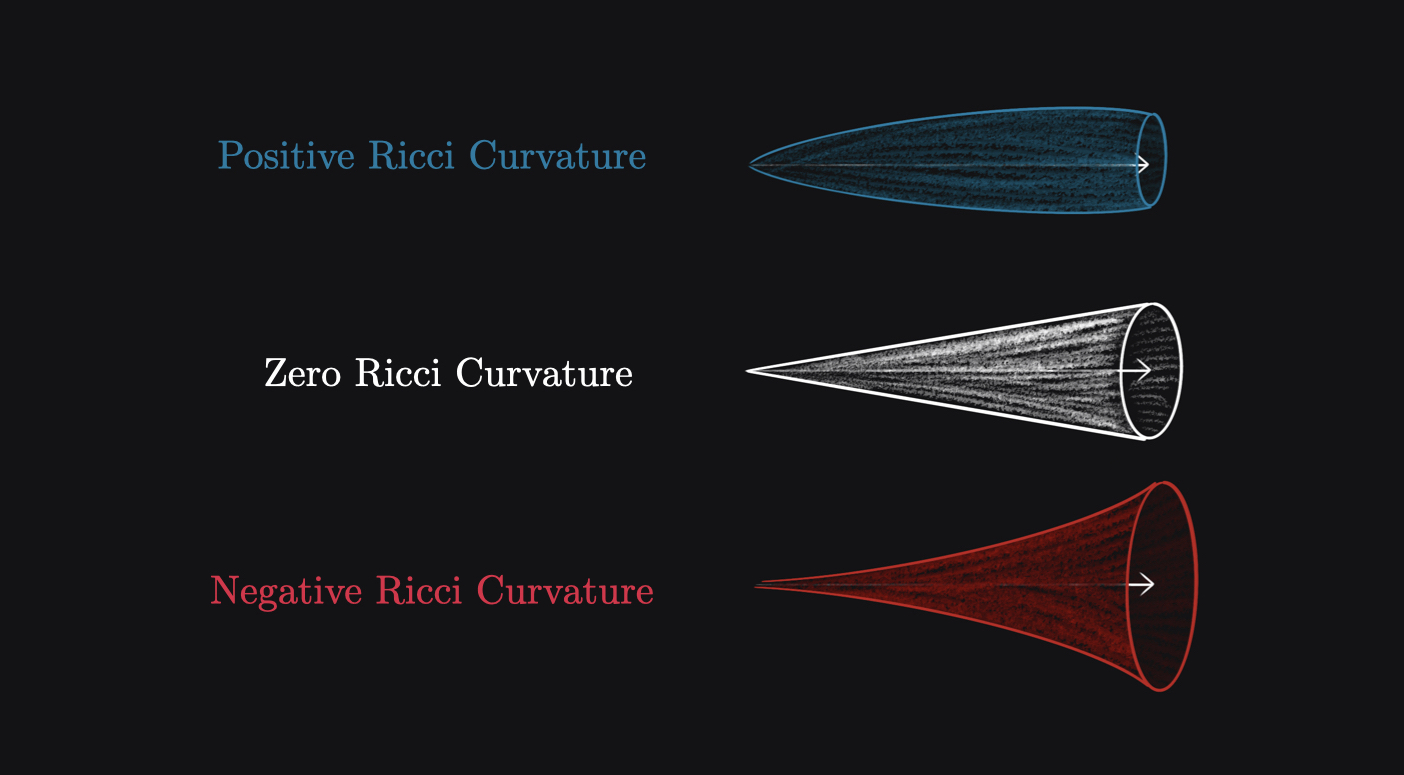}
    \caption{Cones with various Ricci curvatures}
    \label{fig:three cones}
\end{figure}

When the Ricci curvature is positive, the cone will close in on itself whereas in negatively curved space, the cone will open outward (as shown in Figure \ref{fig:three cones}). In this spirit, when we compute the Taylor expansion of $A(\varepsilon),$ we find the following:
\begin{equation} \label{Ricci in Taylor series}
    A(\varepsilon) = \varepsilon^{2(n-1)} \left( D_{n-1} - C_n \Ric(X,X) \varepsilon^2 \right) + O(\varepsilon^{2n+1}).
\end{equation} 

Here $D_{n-1}$ is the volume of a unit disk in Euclidean $(n-1)$-space and $C_n$ is a complicated, but positive, constant which depends on the dimension. The factor of $\varepsilon^{2(n-1)}$ on the right-hand comes from the fact that the cone is both short and narrow, both of which contribute to the area of the base being small.

\subsection*{The Ricci curvature as a geometric Laplacian}

Heuristically, the sectional curvature can be understood as a geometric second derivative. In particular, ignoring the $\varepsilon$ at the front of Equation \eqref{Taylor expansion of L} (which comes from the fact that the associated triangle has very short sides), the sectional curvature appears in the second-order position. With this intuition, comparing Equations \eqref{Laplace as Average} and \eqref{Ricci as average} suggests that the Ricci curvature is a sort of geometric Laplacian. Making this heuristic rigorous requires a bit of Riemannian geometry, because the sectional curvature $K(X,Y)$ does not depend linearly on $X$ and $Y$. 

We won't go into too much detail about how to make this idea precise. The basic idea is that the sectional curvature can be used to construct the \emph{Riemann curvature tensor}, which is comparable to a geometric Hessian.\footnote{It is possible to make the idea of the curvature tensor as a geometric second derivative reasonably precise. In particular, in any set of geodesic normal coordinates, the Riemannian metric satisfies
\[ g_{i j}=\delta_{i j}- \sum_{k,l =1}^n\frac{1}{3} R_{i k j l} x^{k} x^{l}+O\left(|x|^{3}\right) \]
where $R_{i j k l}$ denotes the components of the $(4,0)$ Riemann curvature tensor. The key intuition from this formula is that it strongly resembles a second-order Taylor polynomial for a multivariate function, and the terms $-\frac{1}{3} R_{ijkl}$ correspond to the second-order terms in the expansion.}
 The Ricci curvature is the trace of the Riemann curvature tensor, which gives further credence to its interpretation as a geometric Laplacian.\footnote{There are other ways to make this precise. For instance, in harmonic coordinates, \[\Ric_{ij} = - \frac12 \Delta g_{ij} + \textrm{lower order terms...} \] }

\begin{figure}
    \centering
\includegraphics[width=.9\linewidth]{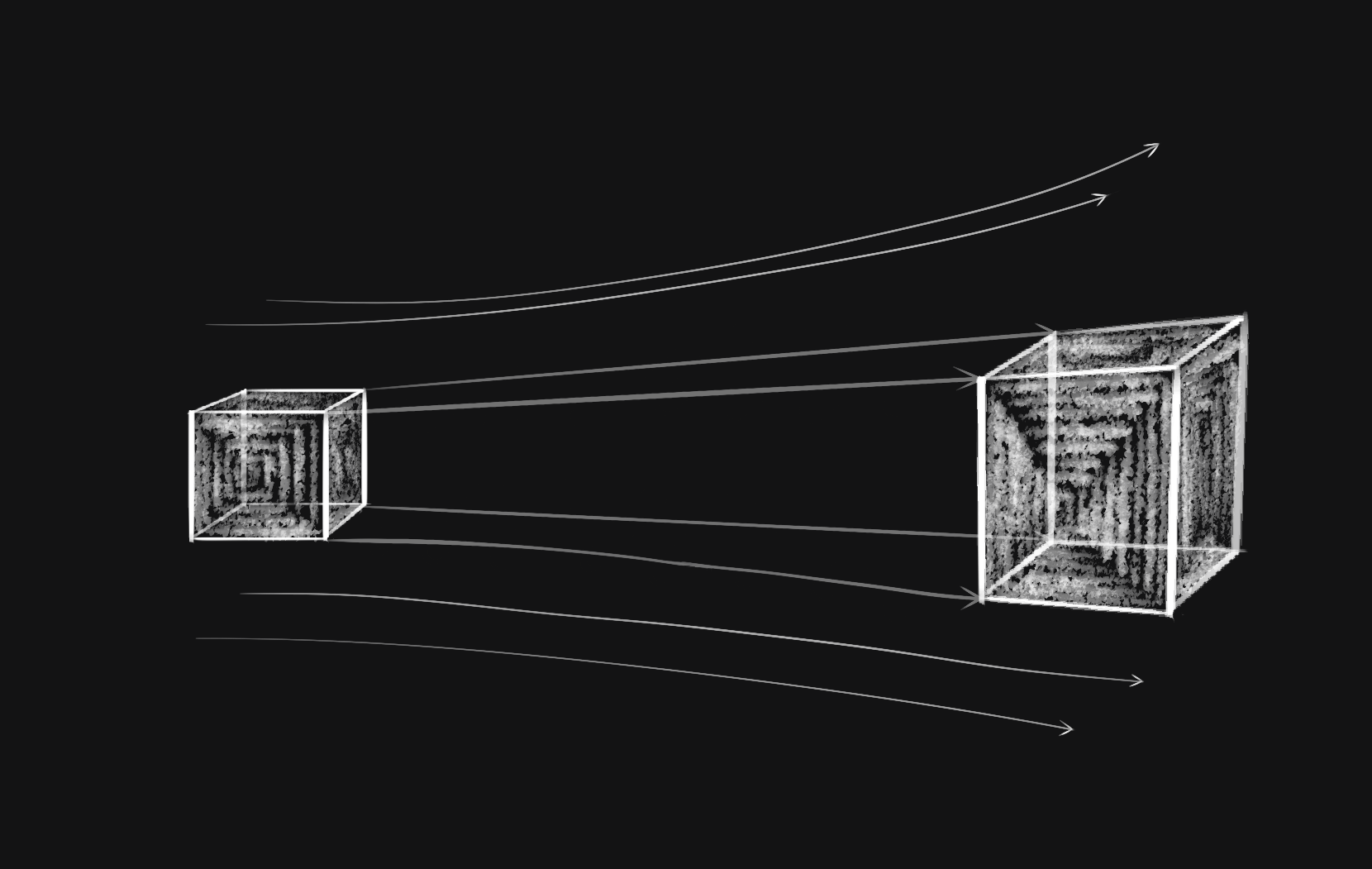}
    \caption[Gradient flow]{The compression of a fluid flow}
    \label{fig:fluid flow}
\end{figure}

There is one other important analogy between the usual Laplacian and the Ricci curvature, which is that they both measure how volumes change. More precisely, given a smooth function $f:\mathbb{R}^3 \to \mathbb{R}$, we can consider the gradient flow
\[ \frac{d}{dt} \gamma(t) = -\nabla f(\gamma), \]
which takes a point and moves in the direction which decreases $f$ the quickest. Doing so, the quantity $\Delta f(\mathbf{x})$ measures how much the volume of a very small cube around $\mathbf{x}$ will change under the flow. In other words, if we consider $\nabla f$ as the current of some fluid, $\Delta f$ measures the compression of the flow. In a similar way, the Ricci curvature determines how volumes of small cubes change as we move from one point to another on a curved space.\footnote{The Weyl tensor is another curvature tensor which is orthogonal to the Ricci curvature and measures the ``tidal forces." In other words, the Weyl tensor determines how the shape of small objects deform when they move along short geodesics whereas the Ricci curvature measures the compression of the gradient flow.} In other words, it measures how volumes are compressed due to the curvature of the space.

\section*{The Ricci Flow}

\label{Ricci Flow section}

Now that we have discussed the heat equation as well as the notion of Ricci curvature, we can finally talk about the Ricci flow, which can be understood as a ``geometric heat equation." To get started, we will provide an informal geometric explanation.

\begin{definition*}
The Ricci flow changes the shape of a space proportional to -2 times the Ricci curvature.
\end{definition*}

\begin{figure}[H]
    \centering
\includegraphics[width=.9\linewidth]{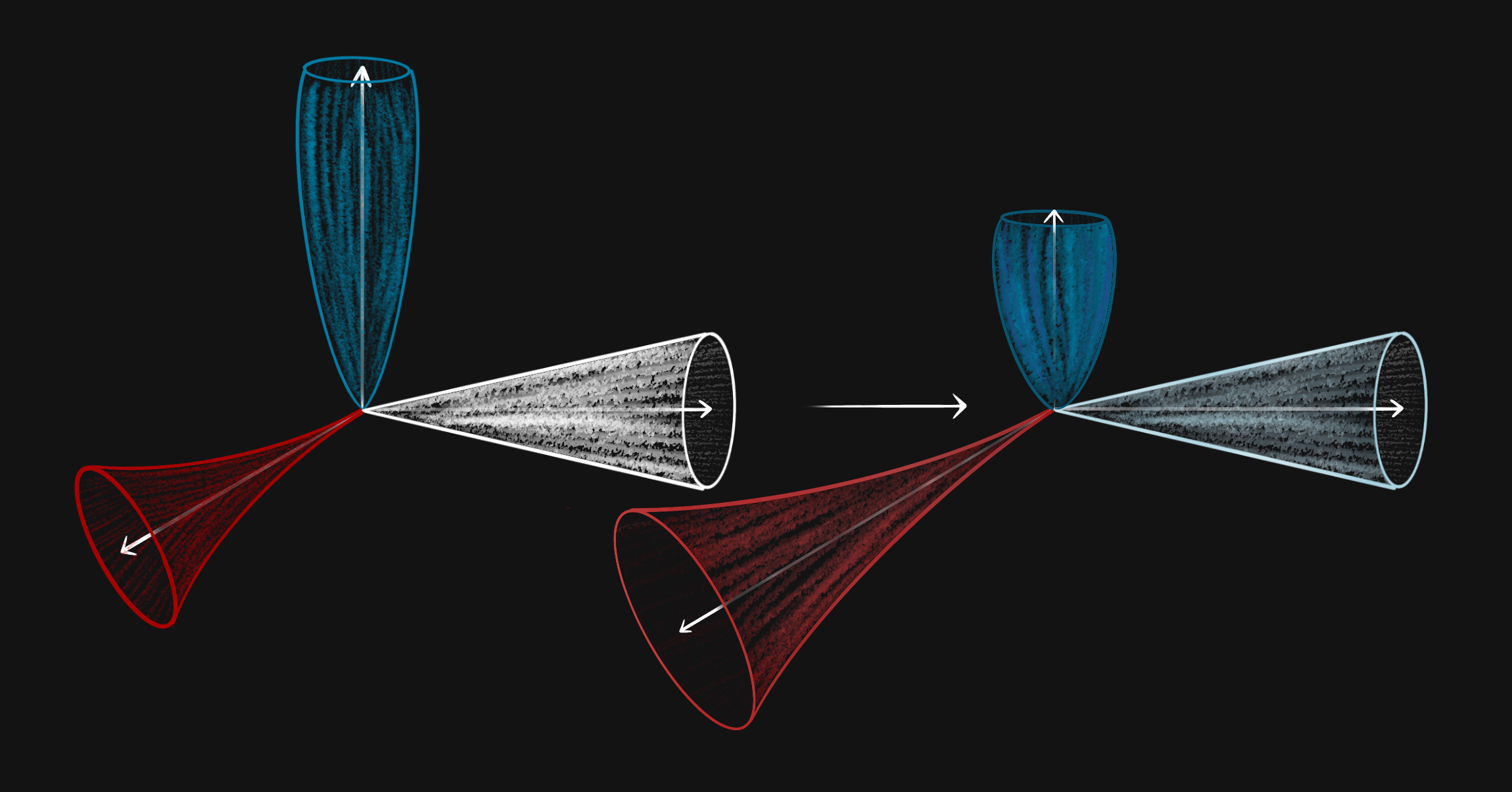}
    \caption[Ricci flow on basis]{Ricci flow deforming a basis which is initially orthonormal\protect\footnotemark}
    \label{fig:Ricciflowbasis}
\end{figure}

\footnotetext{Under Ricci flow, the curvature evolves in a complicated way, which is why the colors of the cones are also changing.}

In other words, directions which have negative Ricci curvature get longer whereas directions with positive Ricci curvature get shorter (as depicted in Figure \ref{fig:Ricciflowbasis}).

In order to define the Ricci flow precisely, we must formalize what the ``shape" of a space is. For surfaces in Euclidean space, it is clear enough what this refers to, but it is much more challenging to define the shape \emph{intrinsically} when it is not lying in some higher dimensional space. To do so, we can use the notion of a Riemannian metric, which is a generalization of the dot product in Euclidean space. In other words, it provides an inner product of two tangent vectors (at the same point) in our space.\footnote{More precisely, given a smooth manifold, a Riemannian metric $g$ is smoothly-varying collection of positive-definite inner products on the tangent space of each point. It is worth noting that a Riemannian metric is not the same as a distance function, which is also known as a \emph{metric}.} There is a fundamental theorem in differential geometry that states that this structure fully determines the geometry of the space, and that it is possible to compute the curvature of the space (as well as all the other metric invariants) from this inner product alone. However, the formulas involved are very complicated, so we will just consider the Riemannian metric $g$ as a geometric object that encodes the ``shape" of our space. Using the metric, the equation for the Ricci flow becomes
\begin{equation} \label{The Ricci flow}
    \frac{\partial}{\partial t} g = -2 \Ric(g).
\end{equation}

It is worth spending some time to make sure this is really a sensible formula, since there are several strange properties.
The first is that our definition of Ricci curvature required that we input a vector (which we wrote, somewhat bizarrely, as two separate entries). In fact, the Riemannian metric $g$ also requires that we input two vectors, so this is actually part of what makes the Ricci flow work.\footnote{Furthermore, both the Ricci curvature and the Riemannian metric are necessarily symmetric, which makes this formula sensible.}

The second objection to considering this equation as a heat equation is that right hand side has a factor of $-2$ whereas the heat equation $\frac{ \partial u}{\partial t} = \Delta u$ does not. This factor appears because the Ricci curvature should be thought of as the \emph{negative} of a geometric Laplacian. In other words, the analytic Laplacian and geometric ``Laplacian" differ by a sign.\footnote{It is important that the coefficient in front of the Ricci curvature is negative, or else the flow will not be defined for forward time (but instead for backwards time).}

\begin{figure}
    \centering
\includegraphics[width=.9\linewidth]{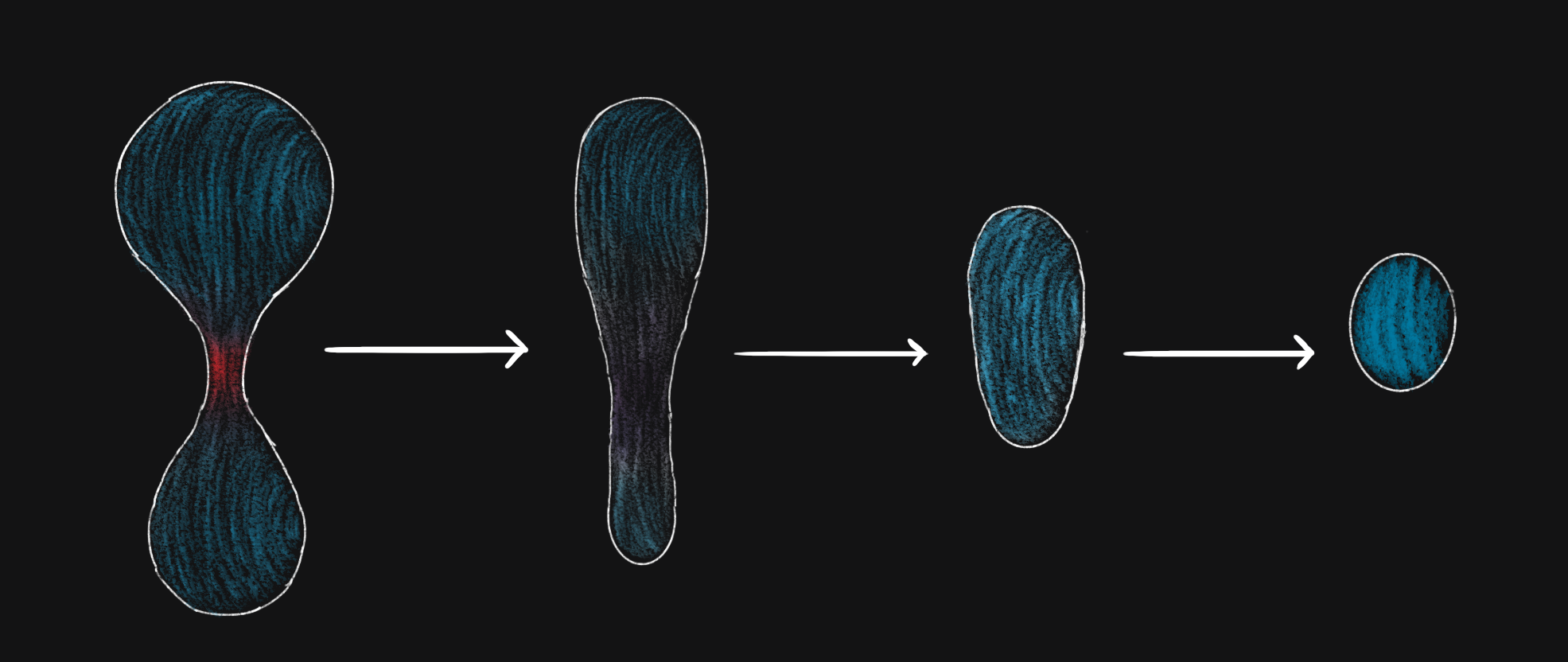}
    \caption{The Ricci flow converging to a round sphere\protect\footnotemark}
    \label{fig:Ricciflowtoroundpoint}
\end{figure}

 With these objections addressed, there is an apparent parallel between this formula and the heat equation $\frac{ \partial u}{\partial t} = \Delta u$. We have tried to justify why it is possible to think of the Ricci curvature tensor as being analogous to a geometric Laplacian, which suggests that the Ricci flow is a heat flow of ``shape." As a demonstration of this fact, it is worthwhile to see an example of the Ricci flow in action, which is depicted in Figure \ref{fig:Ricciflowtoroundpoint}. There are also some animations of the
\href{https://www.youtube.com/watch?v=siAbBsj9XPk}{\textcolor{blue}{flow}}, which are helpful to understand how it behaves ~\cite{Mathifoldvideo}. In both the video and the figure, the surface becomes more spherical as time goes on. In the same way that the heat equation spreads the heat evenly throughout the space, the Ricci flow spreads the curvature evenly throughout the space.

\footnotetext{To be more precise, this flow is converging to a round point, which means that it is shrinking to a point while asymptotically converging to a sphere.}

\section*{Distinctions between the Ricci flow and heat equation}

\label{Distinctions section}

From the analogy between the heat equation and the Ricci flow, we might hope that the Ricci flow will smooth out our space and make it more uniform. In reality, the Ricci flow is more complicated than a heat flow, so this hope is too optimistic.

\begin{figure}
    \centering
\includegraphics[width=.9\linewidth]{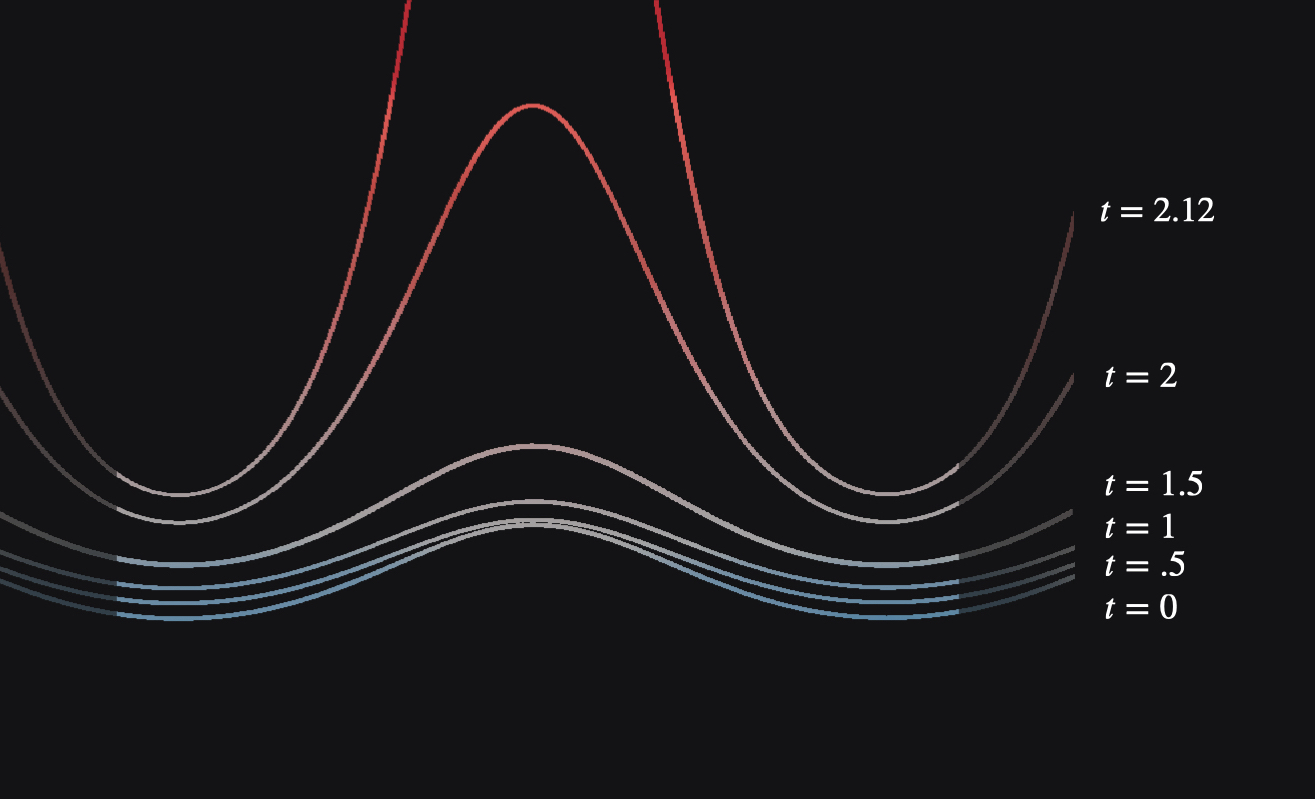}
    \caption{A reaction-diffusion equation going to infinity\protect\footnotemark}
    \label{fig:Quadraticreactiondiffusion}
\end{figure}

The first important distinction between the heat equation and the Ricci flow is that the latter is non-linear (linear combinations of solutions are no longer solutions) because curvature depends in a non-linear way on the metric. Second, it actually behaves more similarly to the reaction-diffusion equation 
\begin{equation} \label{quadratic reaction diffusion}
\frac{ \partial}{\partial t} u = \Delta u + u^2. \end{equation} The first term on the right-hand side behaves as a diffusion term that disperses heat throughout the space whereas the second acts a reaction term that concentrates heat at a point. Reaction-diffusion equations can be thought as a tug-of-war between the diffusion process and the reaction process. If diffusion wins, the solution will smooth itself out much like the normal heat equation. If the reaction term wins out, the heat will become more and more intense and can sometimes even become become infinite in a finite amount of time.\footnote{This occurs with the equation $\frac{ \partial}{\partial t} u = \Delta u + u^2$ if the initial function is positive. For those who have taken a course on PDEs, it is a good exercise to show this using the maximum principle.} An example of this is shown in Figure \ref{fig:Quadraticreactiondiffusion}, where the solution becomes infinitely large after a short amount of time.

\footnotetext{Unlike with the standard heat equation, this equation does not have a closed form solution. As such, we used Fourier analysis to numerically approximate a solution.}

When we compute how the Riemann curvature (denoted $R$) of a space changes along Ricci flow, we find the following equation: \begin{equation} \label{Ricci flow curvature evolution}
     \frac{ \partial}{\partial t} R = \Delta R + R^2  + R^{\sharp}.
 \end{equation}
 It takes some background in Lie algebra to define the terms on the right-hand side of this formula \cite{wilking2013lie}. However, ignoring the final term on the right hand side, there is a notable similarity between this equation and Equation \eqref{quadratic reaction diffusion}.  In particular, sometimes the reaction term wins out and the curvature becomes larger and larger until the shape tears itself apart (or shrinks to a point).\footnote{Once the short-time existence of the Ricci flow has been established, it is not hard to show the flow will exist until the sectional curvature goes to plus or minus infinity. As such, a singularity is generally defined to be a point in space-time where the norm of the sectional curvature tensor goes to infinity.} However, understanding the formation of singularities is difficult because it is not obvious when they occur or what they look like. One of Perelman's key contributions was to classify the possible structures of three-dimensional singularities and to show that a singularity known as the ``cigar soliton" would not develop.  %Intuitively, if you have a space that you want to understand, one natural approach is to consider its Ricci flow and try to show that it converges to a more familiar geometry.

\begin{figure}[H]
    \centering
\includegraphics[width=.9\linewidth]{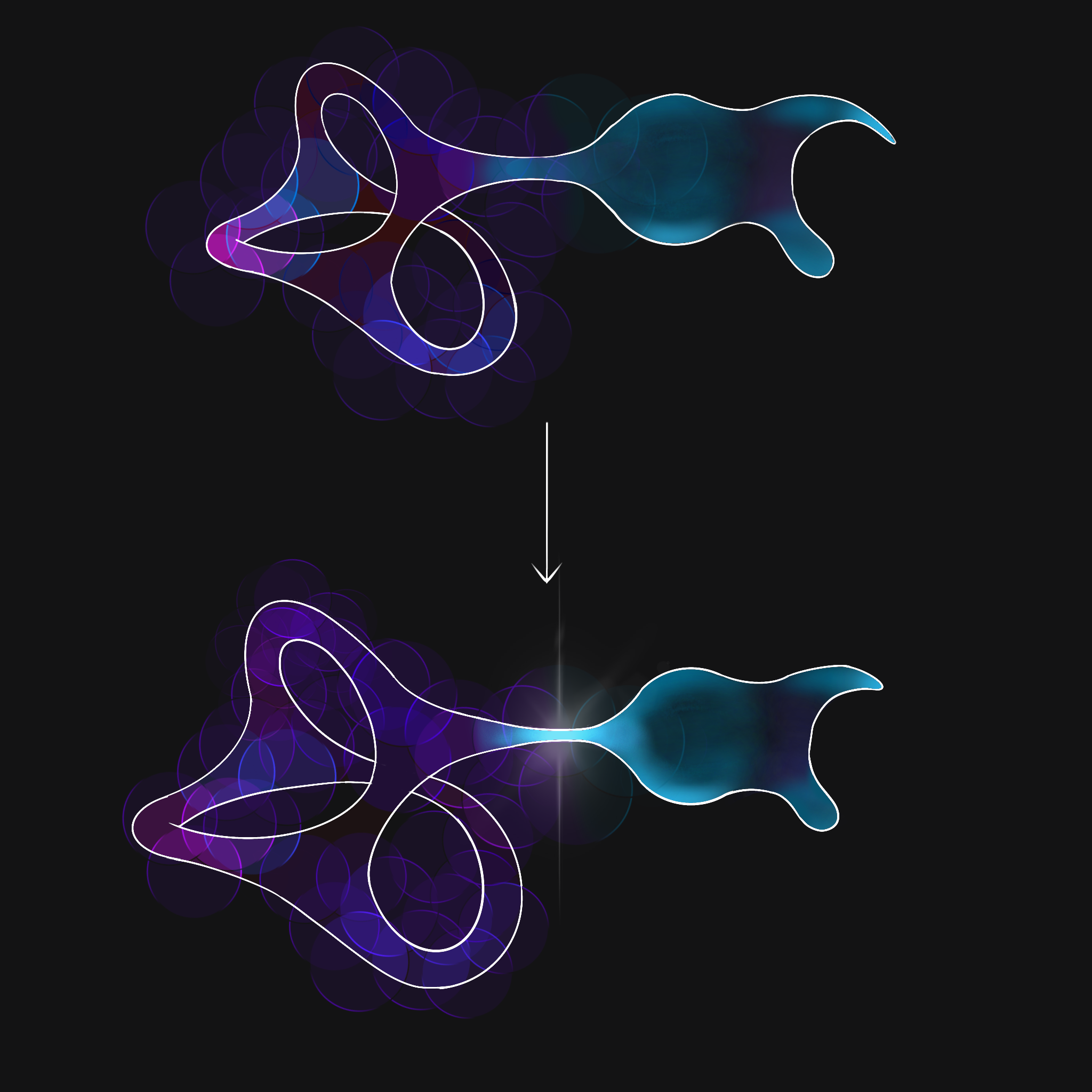}
    \caption{A Ricci flow becoming singular\protect\footnotemark}
    \label{fig:Singular Ricci flow}
\end{figure}

\footnotetext{Neck-like singularities do not occur for the Ricci flow on two-dimensional surfaces, so the space depicted here is three-dimensional. More precisely it is the connected sum of the unit tangent bundle for a hyperbolic surface (which has $\overline{\mathbb{SL}}_2$-geometry) and a deformed $3$-sphere.}

 \chapter*{Recent progress and future directions}
 \addcontentsline{toc}{chapter}{Recent progress and future directions}

 The Ricci flow is an active area of research, both as a tool to prove geometric and topological results and also as a topic of interest in its own right. It would be impossible to give a complete overview of the current state of research, but let us mention a few broad areas of interest.
 
 \subsection*{The long-term behavior and singularity formation of the Ricci flow}
 
The singularity formation and convergence of the Ricci flow is a fascinating and difficult topic. There are many open questions about when singularities occur and what their possible geometries can be. The most famous open problem in this direction is to determine whether the scalar curvature (i.e., the trace of the Ricci curvature) necessarily becomes infinite at a singularity.\footnote{Nata{\v{s}}a \v{S}e\v{s}um  showed that at any singularity, the Ricci curvature becomes infinite~\cite{vsevsum2005curvature}.}
 
There has also been research into the structure of Ricci flows which are allowed to become singular and the singularities that form under ``generic" initial conditions. For details on this line of work, we recommend the recent survey articles by Richard Bamler \cite{bamler2021recent} and Simon Brendle \cite{brendle2022singularity}.

 \subsection*{Geometric classification results}
 
 Apart from the Poincar\'e conjecture, the Ricci flow has been used to prove other geometric classification results. To give a few examples, Hamilton used Ricci flow to understand the geometry of three-dimensional spaces with positive Ricci curvature \cite{hamilton1982three} and four-dimensional spaces with positive curvature operator \cite{hamilton1997four}. In 2007, Brendle and Richard Schoen proved the Differentiable Sphere Theorem using the Ricci flow \cite{brendle2010ricci} and Brendle recently used it to establish a partial classification of higher-dimensional spaces with positive isotropic curvature \cite{brendle2019ricci}.

 There are also several ongoing programs which use Ricci flow (or some related flow) to attack open problems in geometry.
 \begin{enumerate}
     \item The \emph{minimal model program} is a central focus of birational geometry, which is a branch of algebraic geometry involving rational functions. Jian Song and Gang Tian proposed using K\"ahler-Ricci flow (with surgery) to find the minimal models analytically. \cite{song2017kahler}.
     \item The Ricci flow tends to make the curvature of spaces more positive.\footnote{In other words, Ricci flow tends to preserve curvature positivity conditions (for details, see \cite{bamler2019ricci}). As a brief aside, there are a few negative curvature conditions that are preserved as well \cite{khan2020k}.} As such, it is very useful for studying and classifying spaces with positive curvature. For instance, the Ricci flow may be useful for fully classifying spaces which admit metrics of positive isotropic curvature or whose squared-distance has non-negative MTW tensor (see Chapter 12 of \cite{villani2009optimal} for a definition of this condition).
     \item Finally, there are several problems concerning the geometry of four-dimensional spaces which seem amenable to a Ricci flow approach (see Section 5.6 of \cite{bamler2021recent} for some details). 
 \end{enumerate}

 \subsection*{Other geometric flows}
 
 The idea of using a heat-type flow to take geometric space and deform it to some canonical configuration predates the Ricci flow, but this approach has exploded in popularity in the four decades since the Ricci flow was first studied. At present, there are many geometric flows studied in the literature:  mean curvature flow, harmonic map heat flow,  K\"ahler-Ricci flow, Chern-Ricci flow, pluriclosed flow, anomaly flow, Calabi flow, Yamabe flow\ldots
 
These flows play a significant role in differential geometry and have many applications, both within pure mathematics and more broadly. There is still much to be said about the behavior of geometric flows, and it is an active and vibrant area of research.

\begin{comment}

\begin{figure}[H]
    \centering
\includegraphics[width=.9\linewidth]{Images/Torus_shrink_draft.png}
    \caption{Time slices of an Angenent torus shrinking to a point under mean curvature flow}
    \label{fig:Angenenttorus}
\end{figure}
 \end{comment}

\chapter*{Endnotes}
 \addcontentsline{toc}{chapter}{Endnotes}
 
 \section*{Further reading}

Readers who are interested in the Poincar\'e conjecture may enjoy Donal O'Shea's book ``The Poincar\'e Conjecture: In Search of the Shape of the Universe," which is intended for a general audience and gives a detailed history of the problem and its wider role in geometry \cite{oshea2008poincare}.

This paper discussed the uniformization theorem, but did not describe hyperbolic or spherical geometries in detail. For a basic introduction to this subject, I recommend the book ``The Shape of Space" by Jeffrey Weeks. Chapter 6 of Tristan Needham's book ``Visual Complex Analysis" provides an excellent description of the canonical two-dimensional geometries for those who are familiar with complex analysis (or are interested in learning about it) \cite{needham1998visual}. Thurston's ``The Geometry and Topology of Three-Manifolds" is excellent for learning about the three-dimensional geometries \cite{thurston1979geometry}.

For readers who are interested in learning differential geometry, I highly recommend Needham's recent book \cite{needham2021visual} or John Lee's trilogy \cite{lee2010introduction,lee2013smooth, lee2018introduction}. For those conversant with differential geometry who are interested in geometric analysis, Schoen and Shing-Tung Yau's Lectures in Differential Geometry is excellent, although a more challenging read \cite{schoen1994lectures}.

 The discussion of curvature in this paper was strongly influenced by the synthetic theory of curvature bounds. I recommend the following paper of Villani for details on this approach \cite{villani2016synthetic}. There were two reasons I chose this perspective. First, it makes it possible to discuss the geometric meaning of curvature without needing to first define connections, parallel transport, etc. Second, the proof of the Geometrization conjecture uses Alexandrov geometry, so synthetic curvature plays a crucial role in the analysis of Ricci flow. The primary disadvantage of this approach is that it is practically impossible to use for computations, but this was not such an issue for an informal survey.

In order to study the Ricci flow, a good starting point is Peter Topping's manuscript ``Lectures on Ricci flow" \cite{topping2006lectures}. From there, Hamilton's paper on the formation of singularities is well worth reading \cite{hamilton1993formations}. Understanding the full proof of the Poincar\'e conjecture is a massive undertaking (and the full Geometrization conjecture even more so), but there have been several surveys written and I recommend the following \cite{kleiner2008notes, morgan2007ricci}. Furthermore, Danny Calegari recently wrote a chapter on the proof which emphasizes the three-dimensional geometry \cite{calegari2020ricci}.

 \section*{Acknowledgements}

This project started as lecture notes for a talk at the \href{https://math.osu.edu/activities/colloquium/what-is}{\textcolor{blue} {``What is...?" seminar at Ohio State University}} in 2015. The seminar was aimed at advanced high school students at the \href{https://rossprogram.org/}{\textcolor{blue} {Ross program}}, who had taken multivariate calculus and linear algebra, but were not expected to know any differential geometry. I would like to thank the organizers for letting me give a talk about a subject well outside the normal purview of the seminar. 

Finding ways to discuss the Ricci flow without assuming any background in differential geometry was a challenge, and I relied on the help of several people to find explanations which avoided going into detail about Riemannian geometry or PDEs. In particular, thanks to Kori Khan for her helpful suggestions and to Mizan Khan for his help editing the paper. I would also like to thank Frank Nielsen for some corrections.

\bibliography{bibfile}
\bibliographystyle{alpha}

\end{document}